\documentclass[journal ]{new-aiaa}
\usepackage[utf8]{inputenc}
\usepackage{textcomp}

\usepackage{graphicx}
\usepackage{caption}
\usepackage{subcaption}
\usepackage{amsmath}
\usepackage{bm}
\usepackage{algorithm}
\usepackage{algpseudocode}
\usepackage[version=4]{mhchem}
\usepackage{siunitx}
\usepackage{longtable,tabularx}

\newtheorem{thm}{Theorem}[section]

\newtheorem{as}[thm]{Assumption}
\newtheorem{defn}{Definition}[section]

\title{Estimation-Aware Trajectory Optimization\\with Set-Valued Measurement Uncertainties}

\author{Aditya Deole ~\footnote{Ph.D. Candidate,  Department of Aeronautics and Astronautics.} and Mehran Mesbahi ~\footnote{Professor, Department of Aeronautics and Astronautics, AIAA Fellow.}}
\affil{University of Washington, Seattle, Washington, 98195-2400}

\begin{document}

\maketitle
\begin{abstract}
In this paper, an optimization-based framework for generating estimation-aware trajectories is presented. In this setup, measurement (output) uncertainties are state-dependent and set-valued.
Enveloping ellipsoids are employed to characterize state-dependent uncertainties with unknown distributions. The concept of regularity for set-valued output maps is then introduced, facilitating the formulation of the estimation-aware trajectory generation problem. Specifically, it is demonstrated that for output-regular maps, one can utilize a set-valued observability measure that is concave with respect to the finite horizon state trajectories. By maximizing this measure, estimation-aware trajectories can then be synthesized for a broad class of systems. Trajectory planning routines are also examined in this work, by which the observability measure is optimized for systems with locally linearized dynamics. To illustrate the effectiveness of the proposed approach, representative examples in the context of trajectory planning with vision-based estimation are presented. Moreover, the paper presents estimation-aware planning for an uncooperative Target-Rendezvous problem, where an Ego-satellite employs an onboard machine learning (ML)-based estimation module to realize the rendezvous trajectory.
\end{abstract}

\section*{Nomenclature}
{\renewcommand\arraystretch{1.0}
\noindent\begin{longtable*}{@{}l @{\quad=\quad} l@{}}
$\bm x_k \in \mathbb{R}^{n_x},\bm u_k  \in \mathbb{R}^{n_u}, \bm y_k \in \mathbb{R}^{n_y}$  & discrete time state, input, and output triplets at time $k$ \\
$\mathbb{S}^+$ & positive definite matrices (dimension implicit)\\ 
$\mathcal{E}(\bm x,Q)$ & ellipsoid centered at $\bm x$ characterized by $Q \in \mathbb{S}^+$\\ 
$\mathcal{Y}_{\bm x}$ & output set corresponding to state $\bm x$ \\
$\mathcal{X} \subset \mathbb{R}^{n_x},\mathcal{U}\subset \mathbb{R}^{n_u}$ & sets of feasible states and inputs\\
$\Lambda(\mathcal{Y}_{\bm x})$& function defining maximum size of the set $\mathcal{Y}_{\bm x}$\\
$\lVert.\rVert$ & 2-norm \\
$A,B,C$ & state space system matrices \\
$d(.,.),\; d_s(.,.),\; d_{\Gamma}(.,.)$   & metrics for Euclidean, set, and tube distances  \\
$\bm x_{0:T},\; \bm u_{0:T}$ & discrete state and input sequences from time zero to $T$ \\
$\Gamma_{\bm x_0}\bm u_{0:T-1}$ & output tube mapped from initial state $\bm x_0$ with an input sequence $\bm u_{0:T}$\\
$Y_{\bm x_{0:T}}$ & output tube generated for state sequence from time zero to $T$ \\
$D_\mathcal{O}(Y_{{\bm x}_{0:T}}),\;D_\mathcal{O}^{\ell}(Y_{{\bm x}_{0:T}})$ & Degree of observability of an output tube and its lower bound
\end{longtable*}}
\section{INTRODUCTION}

\lettrine{T}RAJECTORY planning often involves solution of optimization problems that efficiently guide a vehicle to a desired final state, while ensuring feasibility with respect to constraints on the vehicle's state and input~\cite{Malyuta22cvx}. This work addresses a trajectory planning problem that accounts for the vehicle's state estimation process—a consideration particularly critical in safety-critical systems where the accuracy of state estimation is essential.

Traditionally, in offline planning, the estimation process is assumed to be independent of state and control, causing the synthesized trajectory to potentially traverse regions with high measurement uncertainty. Meanwhile, in many practical applications, large measurement errors or disturbances adversely affect the quality of the corresponding state estimation; in such scenarios, estimation-aware trajectory design becomes indispensable~\cite{Aware23}.
In order to follow a trajectory, an agent often employs a state estimator, mapping the process measurements to full state estimates, provided that the system is observable~\cite{Kalman60,Gustafsson10filter,Menegaz15filter,JOUIN20162,Taghvaei}. The performance of this estimator is closely tied to the quality of the sensor measurements. If the sensor performance is state-dependent, then visiting states where the sensor has better performance can improve the corresponding state estimation process.
One common approach to design optimal trajectories for improved state estimation is to incorporate a secondary objective into the original trajectory optimization problem~\cite{Alaeddini19Obs,Murali19Oobs}. This augmented objective function must  account for the quality of the estimation process with respect to errors in the measurement. Given knowledge of the measurement error distribution, a suitable metric for capturing the estimation quality can then be developed and embedded in the design objective.
Performance metrics for evaluating the estimation performance include measures on error probabilities~\cite{Blackmore08Prob}, mean squared errors~\cite{Zois17Prob,K07Prob,Zois14Prob}, Fisher information~\cite{Flayac17Fish,Kreutz2013,Fisher_1925}, and entropy and other metrics on the belief state.  
For example, a commonly used construct in settings with a Gaussian measurement noise is the Fisher Information Matrix (FIM)~\cite{Kreutz2013,Fisher_1925}. In a nutshell, FIM quantifies the information content about an unknown state (or parameter) in a measurement. 
In fact in the matrix case, the inverse of FIM provides a lower bound on the expected error covariance, using the Cramer-Rao inequality~\cite{Plastino20,Rao1992}, i.e., $\mathbb{E}\{(\hat{\bm x}-\bm x)(\hat{\bm x}-\bm x)^T\}\geq {\cal I}^{-1}$, where ${\cal I}$ denotes the FIM, and $\bm x$ and $\hat{\bm x}$ denote the state of interest and its estimate, respectively; this lower bound represents the error covariance achieved by an unbiased estimator. 
Various properties of FIM, such its determinant, spectrum, and trace, have been used in different applications~\cite{Singh05}.

The FIM-induced metrics utilized for  quantifying the information content of a measurement scheme are effective for  parameter estimation problems when the underlying state is static. On the other hand, measures such as posterior FIM have been studied to analyze state estimation for linear time-invariant (LTI) systems~\cite{Ober02,Wang2019}. Optimization-based approaches have also been developed to demonstrate the utility of error covariance metrics for trajectory planning problems~\cite{okamoto18,Rafi17}.

Another framework for designing estimation-aware trajectories is via observability-based metrics. Specifically, the Observability Gramian (OG) of a system, mapping initial conditions to an output trajectory sequence, can be used to assess how ``sensitive'' state estimation is to perturbations in output measurements. In this direction, prior studies have demonstrated that properties of the Empirical Gramian (EG) can be leveraged to improve conditioning for the estimation process~\cite{Empirical15, Avnat19Obs}. These works have shown that using EG generated from the model prediction, an optimization problem can be formulated to synthesize a qualitatively more ``observable'' trajectory. These trajectories have improved estimation performance based on the type of metric used on the EG~\cite{Alaeddini19Obs, Alaeddini14Obs, Boyinine22, Glotzbach14}. These studies primarily focus on applications with deterministic output maps with Gaussian noise measurements, that can be propagated forward, multiple roll-outs at a time~\cite{Krener09Unobs,LALL19992598}. 
In addition, probabilistic models using covariance-based optimization have been examined in~\cite{Grebe21,Wilson2014}. Other motion planning approaches use Rapidly Exploring Random Trees to determine estimation-aware trajectories~\cite{Lindemann21}.
And lastly, there are methods proposed in the literature that rely on application-specific correlations between estimation and the underlying state for trajectory design~\cite{Achtelik14}.

In this work, an alternative--yet complementary--perspective to the aforementioned lines of work has been provided, particularly for scenarios where measurements are {\em set-valued} and conventional methods for EG computation are infeasible due to Gaussian noise assumption. However, similar to metrics used EG-based approaches, here, a set distance metric is used to improve the conditioning of the estimation process over a trajectory.
In the area of filtering for set-valued uncertainty, set-membership estimation has been investigated solely for the purpose of filtering, where estimation procedures for set-based outputs are analyzed~\cite{tang20,tang24,Ding}. To the best of our knowledge, studies on estimation-aware planning with set-based output uncertainties have been rather limited, signifying the need for ``filtering'' methods for scenarios where noise distributions are unknown yet bounded. This work specifically focuses on measurement spaces with state-dependent uncertainties that are often non-Gaussian, for example, as they arise in perception-aware planning. 
In this direction, an explicit approach for such planning problems for locally linearized systems is presented. 
Specifically, notions of finite-time observability for nonlinear systems with set-valued state-dependent output uncertainties are first presented.
It is important to point out that in this setup, no specific assumption is made on the distribution for the uncertainties within the set; however, a constraint that the corresponding uncertainty set is compact (i.e., closed and bounded in a finite-dimensional Euclidean space) is enforced. Additionally, it is assumed that a metric (introduced in Section~\ref{sec:Obs_def}), on the size of the uncertainty set has a convexity property with respect to the state vector. This construction, 
 suggests an efficient solution strategy for the corresponding optimization problem. 
In subsequent sections of this paper, justifications for this assumption and its conservatism will be discussed.~\footnote{The size of a compact set can be defined by its diameter or formally as the maximum distance between any two points in the set.}

In this setup, the output map is represented using a bounded set-valued uncertainty over a linear output measurement. 
Although the convexity of the size of this uncertainty with respect to the state is not always guaranteed in real-world applications, smooth variations in the uncertainty and its state dependence are often encountered in practice. In order to address such scenarios, an approximation is proposed for the uncertainty set, using an enveloping convex function that provides a tight upper bound on this set.
Using the convex output map (exact or approximate), a metric is constructed for observability that correlates with the estimation quality. Optimizing this notion of observability is shown to result in better performance for the estimation process. The observability metric is introduced in Section~\ref{sec:Obs_def} as a convex surrogate to manage uncertainty, subsequently incorporated into a generic optimal trajectory design problem, as presented in Section~\ref{sec:case1}.

This paper is also an extension of the previous work reported in~\cite{Deole23}, where a pipeline was proposed to design estimation-aware trajectories using nonlinear model predictive control via distinguishability-based notions. In this work, on the other hand, a framework for establishing guarantees for the existence of estimation-aware trajectories is examined. Specifically, it is shown that a convex optimization problem can be formalized to solve for the desired trajectory when certain conditions on output uncertainty sets are met.

The paper is organized as follows: In Section~\ref{sec:model}, the system and assumptions on properties of the set-valued output uncertainty map are presented. In Section~\ref{sec:Obs_def}, observability conditions are introduced and the primary criterion for finite horizon weak observability for set-valued maps are presented. The notion of \textit{degree of observability}, serving as an evaluation metric for the planned trajectory, is then examined. Section~\ref{sec:Obs_thm} presents the main results of this work, including a discussion pertaining to the convexity of the observability-based metrics. Next, in Section~\ref{sec:traj}, the optimization problem is formulated to compute the desired deviations, relative to a nominal trajectory, to maximize the degree of observability. Finally, in Section~\ref{sec:case1}, two scenarios are presented to demonstrate the applicability of the theoretical developments. In the first example, it is shown that for a predefined trajectory (to achieve a certain task), trajectory deviations can be characterized to improve the state estimation
when the output maps are \textit{observability-regular}.
A nonlinear double-integrator Dubins car model is used to showcase this scenario.
For the second scenario
presented in Section~\ref{sec:sat}, an estimation-aware trajectory is designed alongside the task optimization problem. In this example, an Ego-Target Rendezvous problem is considered to showcase an iterative optimization approach built upon the results presented in~Section~\ref{sec:traj} for generating estimation-aware trajectories. Further implementation details about the sequential programming approaches are presented in the Appendix. An implementation of a basic trust region algorithm is presented in Section~\ref{BTR} and details of the approach used to solve the satellite rendezvous problem are presented in Section~\ref{sec:scvx_imp}.  

\section{Set-valued uncertainty measurement model}\label{sec:model}

This section formalizes an observability-based metric that is subsequently used to design estimation-aware trajectories. In particular, this metric is developed such that maximizing its {\em lower bound} improves observability of the state, hence improving the measurement-driven estimation process. Subsequently, it is shown that this metric can be augmented to the objective of a baseline trajectory generation algorithm.

The model that is considered for the trajectory planning problem is built upon the discrete-time nonlinear system of the form,
\begin{align}
    {\bm x}_{t+1} &= f({\bm x}_t,{\bm u}_t), \label{Eq:sys1} \\
    {\bm y}_t &= h({\bm x}_t), \label{Eq:sys2}\\
    {\bm y}_t' &= h({\bm x}_t) +  {\bm \varepsilon}_t, \label{Eq:sys3}
\end{align}
with $n_x,n_u$ and $n_y$ as the state, input and output dimensions;
in this case, $\bm x_t \in \mathcal{X} \subset \mathbb{R}^{n_x}$ and $u_t \in \mathcal{U} \subset \mathbb{R}^{n_u}$ are feasible state and control inputs for $t=0,1,2, \ldots,T$ and the term $\varepsilon_t$
are set-valued, i.e., $\varepsilon_t \in {\cal E}$, as opposed to say, point-wise deterministic measurements or sampled from a known probability distribution. In particular, the output ${\bm y}_t'$
is sampled uniformly from a compact set centered around the nominal output $ {\bm y}_t  = h({\bm x}_t)$. The functions $f:\mathbb{R}^{n_x}\times \mathbb{R}^{n_u} \rightarrow \mathbb{R}^{n_x}$ and $h:\mathbb{R}^{n_x} \rightarrow \mathbb{R}^{n_y}$ are assumed to be differentiable.

The set of possible outputs with respect to the measurement uncertainty set ${\cal E}$ is defined as
$\mathcal{Y}_{\bm x}=h(\bm x) \oplus {\cal E}$,
where ``$\oplus$'' denotes the set addition;\footnote{The set generated by adding $h(\bm x)$ to all elements of ${\cal E}$.} 
${\cal E}$ encodes the uncertainty in the observations due to noise, disturbances, or other unknown parameters effecting the output.
In this work, the measurement uncertainty is characterized by maximal ellipsoids for which axial radii can be computed. The size of the uncertainty set is then defined as the maximum ellipsoidal radius. A ``variation'' 
in this notion of uncertainty is then utilized in the subsequent optimization formulation of estimation-aware planning. The maximal ellipsoid will bound the set-valued uncertainty observed in the model by a state-dependent set centered, around the nominal measurement 
$\bm y=h(\bm x)$.
A general state-dependent ellipsoidal set centered at $\bm x_0$ is defined as, 
\begin{align}
\mathcal{E}({\bm x_0},Q(\bm x_0)) = \left\{ \bm x \;|\; (\bm x - \bm x_0)^\top \, Q(\bm x_0)^{-1} \, (\bm x - \bm x_0) \leq 1\right\}, \label{Eq:def_ellipse}
\end{align}
where $Q(\bm x_0) \in \mathbb{S}^+$ and is dependent on ${\bm x}_0$; occasionally, the short-hand notation $\mathcal{E}({\bm x}_0)$ for a generic ellipsoid centered at ${\bm x}_0$. For $\bm x_0 \in \mathcal{X}$ is adopted. 
The point-to-set map is defined as, 
\begin{align}
    \mathcal{Y}_{\bm x_0} := \mathcal{E}({\bm y_0},Q(\bm x_0)), \label{Eq:ellipse}
\end{align}
where $\mathcal{Y}_{\bm x_0}$ denotes the ellipsoidal output uncertainty around 
$\bm y_0 = h(\bm x_0)$; note that the ellipsoidal output is a dependent on the corresponding state. Subsequently, $\bm y$ is used to denote the deterministic output corresponding to the state $\bm x$, and $\mathcal{Y}_{\bm x}$ as the corresponding output uncertainty set. 
The map from the state $\bm x \in \mathcal{X}$ to the set $\mathcal{Y}_{\bm x}$ defined in Eq. (\ref{Eq:ellipse}), is denoted simply by $\mathcal{Y}$ for brevity
Ellipsoidal sets are often used for representing uncertainty sets~\cite{Ellipsoidal}; in this work, these ellipsoids are used to represent and approximate uncertain set-valued output measurements in this work.

The key problem addressed in this work is the selection of the input sequence such that the corresponding state sequence remains ``distinguishable'' throughout its trajectory. Notions of distinguishability for characterizing the observability of a system have been examined by Hermann and Krener for a general class of deterministic systems~\cite{Hermann77Obs}. In this work, these notions are modified to reason about the observability of uncertain set-valued output sequences. Note that in this setting, the map of state to output uncertainty ellipsoids is assumed to be known. Otherwise, its action has to be validated in a simulation environment. This assumption is analogous to the case where the error covariance of a sensor is assumed to be given in the filter design.
For the optimization approach to estimation-aware planning, a metric is required to quantify the size of the uncertainty set. Here, a conservative bound for comparing ellipsoidal uncertainty sets is adopted, based on their radii. Specifically, the maximum eigenvalue of the positive definite matrix defining the ellipsoid will be used as this metric.

Note that for an ellipsoid centered at $\bm x$, the distance from $\bm x$ to any other element of the ellipsoid is bounded by the largest eigenvalue of the positive definite matrix defining the ellipsoid. Hence,
\begin{align}
        \sup\limits_{\bm y \in \mathcal{Y}_{\bm x}} d(\bm x,\bm y)   = \lambda_{\max}(Q_{\bm x}) ; \label{Eq:maxellipse}
    \end{align}     
the notation $\Lambda(\mathcal{Y}_{\bm x})$ denotes the largest eigenvalue of the ellipsoidal uncertainty corresponding to state $\bm x$; see Fig.~\ref{fig:map1}. In this figure, we show the map $\mathcal{Y}$ transports $\bm x$ to its uncertainty set $\mathcal{Y}_{\bm x}$, shown in blue, and the function $\Lambda(\mathcal{Y}_{\bm x})$ covers the uncertainty with its maximum size shown in gray. 
%

Given this measure of ellipsoidal uncertainty for a state,
its variation as a function of the state can be quantified.
A key assumption in the subsequent analysis is as follows.

\begin{as}\label{as:eigs}
    Let the largest eigenvalue of the output uncertainty ellipsoid corresponding to state $\bm x$ be denoted by $\Lambda(\mathcal{Y}_{\bm x})$. The map $\Lambda$ is convex and uniformly bounded with respect to $\bm x \in \mathcal{X}$. 
\end{as}
 
Convexity of the above measure of uncertainty implies Lipschitzness of $\Lambda(Q_{\bm x})$~\cite{Roberts74}; as such, for any feasible $\bar{\bm x},\bm x \in \mathcal{X}$ one has,
\begin{align}
|\Lambda(\mathcal{Y}_{\bm x})-\Lambda(\mathcal{Y}_{\bar{\bm x}})| \leq L(\bar{\bm x},r) \, \lVert \bm x-\bar{\bm x}\rVert. \label{Eq:lips}
\end{align}
For the purpose of subsequent discussion, the notation $\bar{\bm x}$ in Eq. (\ref{Eq:lips}) signifies a state on the nominal trajectory; the state
$\bm x$ on the other hand, belongs to an open neighborhood $U$ of $\bar{\bm x}$.
In Eq. (\ref{Eq:lips}), $L(\bar{\bm x},r)$ is the ``local'' Lipschitz parameter. This parameter can be chosen as,
\begin{align}
L(\bar{\bm x},r) = \frac{2M(\bar{\bm x},r)}{r}, \quad\text{where} \quad M(\bar{\bm x},r) = \max\limits_{\bm x \in B_{r}(\bar{\bm x})}\left\{\Lambda(\bm x)\right\} ;\label{eq:lips2}
\end{align}

in Eq. (\ref{eq:lips2}), $M$ represents the local bound on $\Lambda(\bar{\bm x})$ in the neighborhood $B_{r}(\bar{\bm x})$ of $\bm x$. One can observe that $L(\bar{\bm x},r)$ is convex with respect to $\bar{\bm x}$.
Note that if the analytical expression for $\Lambda(\bar{\bm x})$ is given, then $L(\bar{\bm x},r)$ can be obtained via (\ref{eq:lips2}) for a given radius $r$.
 For certain applications where a perfect sensor model is not known but 
Assumption~\ref{as:eigs} holds, both $L(\bar{\bm x},r)$ and  $\Lambda(\bar{\bm x})$ can be characterized via a simulation oracle.

Equation (\ref{Eq:lips}), gives the inequality,
\begin{align}
|\Lambda(\mathcal{Y}_{\bm x})| \leq |\Lambda(\mathcal{Y}_{\bar{\bm x}})| + L(\bar{\bm x},r)\lVert \bm x-\bar{\bm x}\rVert;\label{Eq:Lbound}
\end{align}
picking $r = \lVert \bar{\bm x} - \bm x\rVert$, one can obtain a viable yet conservative bound on $|\Lambda(\mathcal{Y}_{\bm x})|$. For constructing the set-valued observability metric, set $r = \delta\bm x$, for the maximum deviation $\delta {\bm x}$ from the trajectory for which the metric is being computed. Assumption~\ref{as:eigs} and inequality Eq. (\ref{Eq:Lbound}) facilitate ``control'' over the size of the variation
in output ellipsoidal uncertainty. \\

\subsection*{A note on convexity of the output map}

The generation of estimation-aware trajectories can be mathematically formalized when one can quantify how certain states improve the state estimation process relative to others. Conversely, if the uncertainty set is uniform over all states, estimation -aware planning becomes unnecessary, as it does not improve observability. For cases with a known state dependency on the state, the aim is to design an optimization problem that produces a unique exploration trajectory. In Section~\ref{sec:Obs_def}, it is shown that Assumption \ref{as:eigs} plays a crucial role in ensuring a unique solution to this planning problem.

In some applications, state dependence may not exhibit a convexity property but instead display a local monotonicity that aids exploration. For example, quasi-convexity of state-dependent uncertainties become relevant in such cases; this is common in vision-based sensing, for example, where identifiers or features are most visible at specific states, with sensing accuracy diminishing monotonically as one moves away from these optimal visibility states.
To address such scenarios, the size of the uncertainty set is approximated using an enveloping convex function.
This function is defined as
\begin{align}
  \hat{\Lambda}(\bm x) = \inf \{g({\bm x}) \,| \, g \text{ is convex and }g({\bm x})\geq \Lambda({\bm x})\},  \label{eq:lambdaapprox}
\end{align}
for ${\bm x}$ in a domain of interest.
The selection of candidate convex function $g$ depends on specific applications.~\footnote{In the context of this work, the enveloping function $g$ is required be convex. The convex metric derived from this function will be used as the objective in an optimization problem. For faster convergence rates with say gradient descent, one can select $g$ to be strongly convex.} Note that this approach results in an over-approximation of the upper bound. In subsequent analysis, notation $\Lambda$ refers to the envelope $\hat{\Lambda}$.

Using the enveloping function, the observability metric is computed using a separation function. This separation function is a distance metric between sets, defined to be positive when ``distinguished'' trajectories are generated from neighboring initial conditions. 
In practice, designing the enveloping function from the data is referred to as the validation procedure. For the case study on the Ego-Target problem, this validation procedure is detailed in Sections~\ref{app:sim} and \ref{sec:experiment}.   

\subsection*{Set Distance metric}
For the upcoming observability analysis, one has to define 
a notion of separation distance between two sets.
\begin{defn}\label{def:setdist}
(\textit{Set Distance}) The distance between sets $\mathcal{A}$ and $\mathcal{B}$ in $\mathcal{X}$ 
is defined as,
\begin{align}
    d_s(\mathcal{A},\mathcal{B}) = \inf\left\{ d(\bm v,\bm w)\; | \; \bm v\in \mathcal{A},  \bm w \in \mathcal{B} \right\} , \label{Eq:setdistdef}
\end{align}
\end{defn}
where $d(\bm v,\bm w)$ for a pair of elements in ${\bm v}, {\bm w} \in {\cal X}$ is their Euclidean distance.

Note that Eq. (\ref{Eq:setdistdef}) defines a \textit{pseudo-metric} on sets in $\mathcal{X},$ as  $d_s(\mathcal{A},\mathcal{B}) = 0$ does not imply that $\mathcal{A}$ and $\mathcal{B}$ are the same set; instead, this distance
is zero when $\mathcal{A}\cap\mathcal{B} \neq \emptyset$. In this paper, for simplicity, the pseudo-metric referred as a metric, only on separated sets. It is implied that metric properties of the distance function Eq. (\ref{Eq:setdistdef}) holds when sets are separated and that $d_s(\mathcal{A},\mathcal{B}) = 0$, also written as $\mathcal{A} = \mathcal{B}$, implies their \textit{non-separability} instead of them being identical. The same interpretation holds for the tube distance ``metric'' defined subsequently. As such, in this analysis, maximizing set separation is only applied for separated sets.

The set distance Eq. (\ref{Eq:setdistdef}) is used to quantify separation between the uncertainty ellipsoids; in particular, a nonnegative lower bound for this distance uses Eq. (\ref{Eq:setdistdef}). In this direction, consider the sets $\mathcal{E}(\bm x_1)$ and $\mathcal{E}(\bm x_2)$, with $\bm v \in \mathcal{E}(\bm x_1)$ and $\bm w \in \mathcal{E}(\bm x_2)$.
Since, 
\begin{align}
     d(\bm x_2,\bm v) &\leq d(\bm v,\bm w) + d(\bm w,\bm x_2), \label{Eq:dist1}\\
     d(\bm x_1,\bm x_2) &\leq d(\bm x_1,\bm v) + d(\bm x_2,\bm v), \label{Eq:dist2}
\end{align}

by substituting Eq. (\ref{Eq:dist1}) in Eq. (\ref{Eq:dist2}) we obtain,
\begin{align}
    d(\bm v,\bm w) &\geq d(\bm x_1,\bm x_2) - d(\bm x_1,\bm v) - d(\bm w,\bm x_2)>0, \, \label{Eq:dist3}
\end{align}
where the last inequality for positivity implies that these sets are separated. Taking the infimum in Eq. (\ref{Eq:dist3}), it now follows that,
 \begin{align}
    \inf\limits_{\bm v,\bm w} \; d(\bm v,\bm w) &\geq \inf\limits_{\bm v,\bm w} \left\{d(\bm x_1,\bm x_2) - d(\bm x_1,\bm v) - d(\bm w,\bm x_2)\right\} \nonumber \\
    &\geq  d(\bm x_1,\bm x_2) -\sup\limits_{\bm v}\left\{d(\bm x_1,\bm v)\right\} - \sup\limits_{\bm w} \left\{d(\bm w,\bm x_2) \right\} . \label{Eq:dist4}
\end{align}
 A simple lower bound on this set distance is given by,
\begin{align}
    d_s(\mathcal{E}(\bm x_1),\mathcal{E}(\bm x_2)) &\geq \lVert \bm x_1 - \bm x_2 \rVert_2 - \lvert \Lambda(Q(\bm x_1)) - \Lambda(Q(\bm x_2))  \rvert . \label{Eq:setdistlb}
\end{align}

In order to ensure that two ellipsoidal sets with centers $\bm x_1,\bm x_2,$ are \textit{separable}, we will enforce the inequality,
\begin{align}
    d(\bm x_1,\bm x_2) \geq \Lambda(Q(\bm x_1)) + \Lambda(Q(\bm x_2));\label{eq:separable}
\end{align}
this will ensure that the corresponding set distance Eq. (\ref{Eq:setdistlb})
is in fact a pseudo-metric.

\section{Observability with set-valued output uncertainty}\label{sec:Obs_def}

To ``control'' the quality of the estimation process along a planned state sequence $x_{0:T}$, a set-valued observability measure is defined in this section. For a LTI system with a Gaussian noise model, it is known that estimation can be improved by directly optimizing properties of the OG~\cite{Alaeddini14Obs}. For set-valued measurements, an analogous setup is proposed in this work. In this direction, the notion of \textit{distinguishability} is first examined~\cite{Hermann77Obs}.

Consider the system Eq. (\ref{Eq:sys1})-(\ref{Eq:sys2}) with a given input sequence $\bm u_{0:T-1}$ and an initial condition $\bm x_0$. Here, the output sequence $\bm y_{0:T}$ is defined as the image of the map $\Sigma$ parameterized with a fixed input sequence for any initial condition. Thereby, define the initial condition to output sequence mapping as $\Sigma_{{\bm x}_{0}}(\bm u_{0:T}) = \bm y_{0:T}$.
We say that state $\bar{\bm x}_0$ is \textit{indistinguishable} from another state ${\bm x}_0 \in \mathcal{X}$ if for every admissible input sequence, we have $\Sigma_{\bar{\bm x}_0}(\bm u_{0:T}) = \Sigma_{\bm x_0}(\bm u_{0:T})$.
The system is said to be observable at $\bar{\bm x}_0$, if $\bar{\bm x}_0$ is distinguishable. \textit{Distinguishability} ensures that for 
some input sequence and distinct initial conditions ${\bm x}_0, \bar{\bm x}_0 \in \mathcal{X}$, $\Sigma_{\bar{\bm x}_0}(\bm u_{0:T}) \neq \Sigma_{\bm x_0}(\bm u_{0:T})$, i.e.,
$ d_s(\Sigma_{\bar{{\bm x}}_0}({\bm u}_{0:T}),\Sigma_{{\bm x}_0}({\bm u}_{0:T}) > 0$.
\begin{figure}[t]
\centering
\subcaptionbox{State to output uncertainty mapping \label{fig:map1}}
{\includegraphics[width = 0.3\linewidth]{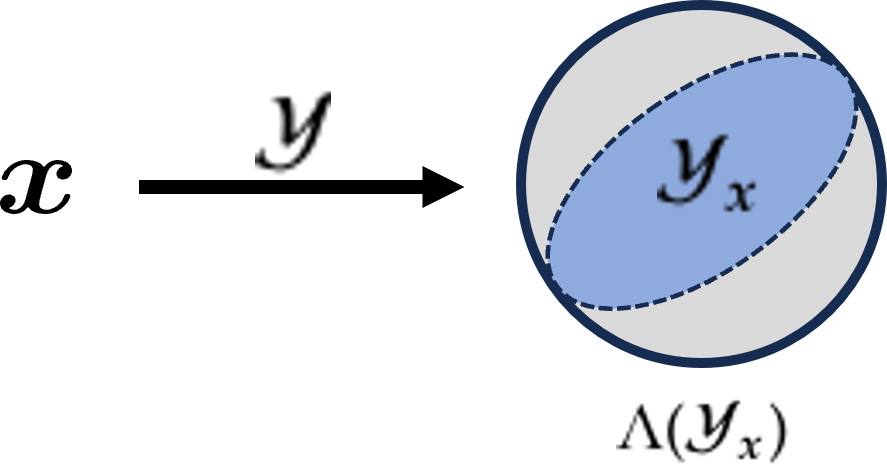}}
\subcaptionbox{Initial condition to output tube mapping\label{fig:map2}}
{\includegraphics[width=\linewidth]{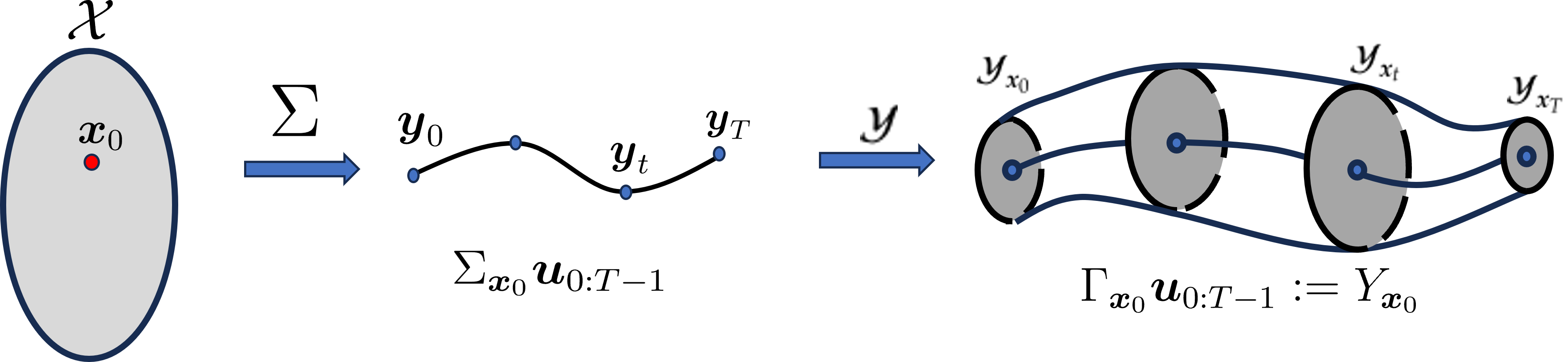}}
\caption{Representations of static and dynamics state to output maps for a given input sequence.}\label{fig:maps}
\end{figure}

The notion of set separation is now extended to set-valued {\em output} sequences that will be referred to as ``tubes.''
Such output sequences are characterized via the composition of the point-to-set map $\mathcal{Y}$ with $\Sigma$ denoted by $\mathcal{Y}\circ\Sigma:=\Gamma$. 
As such, $\Gamma$ maps the initial condition $\bm x_0$ to the output tube $\{\mathcal{Y}_{x_0},\mathcal{Y}_{x_1},\ldots,\mathcal{Y}_{x_T}\}$ for an input sequence $\bm u_{0:T-1}$; in this paper, the tube is a time-series of uncertainty sets denoted by the variable $Y_{\bm x_{0:T}} = \{\mathcal{Y}_{x_0},\mathcal{Y}_{x_1},\ldots,\mathcal{Y}_{x_T}\}$. 
For a fixed initial condition, the mapping from initial condition to output tube is illustrated as follows:
\begin{align}
    \bm x_0 \; \underbrace{\xrightarrow{\Sigma(\bm u_{0:T-1})} \; \Sigma_{\bm x_0} \, \bm u_{0:T-1} \quad \xrightarrow{\quad\mathcal{Y}\quad}}_{\Gamma=\mathcal{Y}\circ\Sigma} \quad \Gamma_{\bm x_0}\bm (u_{0:T-1})\; := \;Y_{\bm x_{0:T}}\label{Eq:map_def}
\end{align}
This mapping scheme is illustrated in Fig.~\ref{fig:map2}.
The first map (depicted on the left panel) generates an output sequence $\Sigma_{\bm x_0}\bm u_{0:T-1}$ for some $\bm x_0 \in \mathcal{X}$. Then this sequence is converted to the tube sequence shown in Fig.~\ref{fig:map2} (right panel), with ellipsoidal sets generated for each point $\bm y_t$ with the size dependent on ${\bm x}_t$. 

Now the notion of observability is defined such that, $\Gamma$ is a one-to-one mapping from ${\bar{\bm x}}_0$ to the output tube of length $T$;
for brevity, $\bm u_{0:T-1}$ will be dropped from the $\Gamma$ map.

 \begin{defn}\label{tube_dist}
The distance metric $d_{\Gamma}(.,.)$ on output tubes $Y_{\bar{\bm x}_{0:T}}$ and $Y_{\bm x_{0:T}}$ over the time interval $[0,T]$ ($T>0$), generated by initial conditions $\bar{{\bm x}}_{0}$ and ${\bm x}_{0}$ for an input sequence, is defined as, 
\begin{align}
d_{\Gamma}\left(Y_{{\bar{\bm x}}_{0:T}},Y_{\bm x_{0:T}}\right) = \sum\limits_{t=0}^{T} d_s\left(\mathcal{Y}_{{\bar{\bm x}}_t},\mathcal{Y}_{{\bm x}_t}\right), 
\end{align}
where the set-distance is given by Definition~\ref{def:setdist}.~\footnote{Note that $d_{\Gamma}$ is a pseudo-metric as $d_{\Gamma}\left(Y_{{\bar{\bm x}}_{0:T}},Y_{\bm x_{0:T}}\right) = 0$ does not imply that the tubes $Y_{{\bar{\bm x}}_{0:T}}$ and $Y_{\bm x_{0:T}}$ are identical. In fact, the metric is only applicable for tubes where $Y_{{\bar{\bm x}}_{0:T}}\bigcap Y_{\bm x_{0:T}} \neq \emptyset \implies  d_{\gamma}(.,.)>0$.}
\end{defn}
The above notion of tube distance is analogous to that of sets; as such, 
it is also a pseudo-metric since $d_{\Gamma}\left(Y_{{\bar{\bm x}}_{0:T}},Y_{\bm x_{0:T}}\right)=0 $ does not imply that tubes $Y_{{\bar{\bm x}}_{0:T}}$ and $Y_{\bm x_{0:T}}$ are identical. 
The function $d_{\Gamma}$ is now used to define the following notion of observability.
\begin{defn}\label{Obs}
The system with a set-valued output is said to be \textit{finite horizon weakly observable} at ${\bar{\bm x}}_0$ over horizon $T>0$, if for all $\bm x$ in an open neighborhood $U$ of ${\bar{\bm x}}_0$, and an arbitrary input sequence, one has,
\begin{align}
d_{\Gamma}\left(Y_{{\bar{\bm x}}_{0:T}},Y_{\bm x_{0:T}}\right) = \sum\limits_{t=0}^{T} d_s\left(\mathcal{Y}_{{\bar{\bm x}}_t},\mathcal{Y}_{{\bm x}_t}\right) \;> 0, \quad \mbox{implies that} \quad d({\bar{\bm x}}_0,\bm x_0)>0; \label{Eq:dint}  
\end{align}
as such, $d_{\Gamma}\left(Y_{{\bar{\bm x}}_{0:T}},Y_{\bm x_{0:T}}\right)>0$ implies that $\bar{\bm x}_0$ and $\bm x_0$ are distinguishable.  
\end{defn}
For tubes to be \textit{separable}, at least one output-uncertainty pair ($\mathcal{Y}_{{\bar{\bm x}}_t},\mathcal{Y}_{{\bm x}_t}$), generated using the same input sequence, should be \textit{separable} for some $t$.
In this case, Eq. (\ref{Eq:dint}) holds and trajectories ($Y_{{\bar{\bm x}}_{0:T}},Y_{\bm x_{0:T}}$) are \textit{distinct}, necessary for the validity of set-valued observability. 

\begin{figure}[t]
    \centering
    \smallskip
    \smallskip
    \smallskip
    \includegraphics[width=\linewidth]{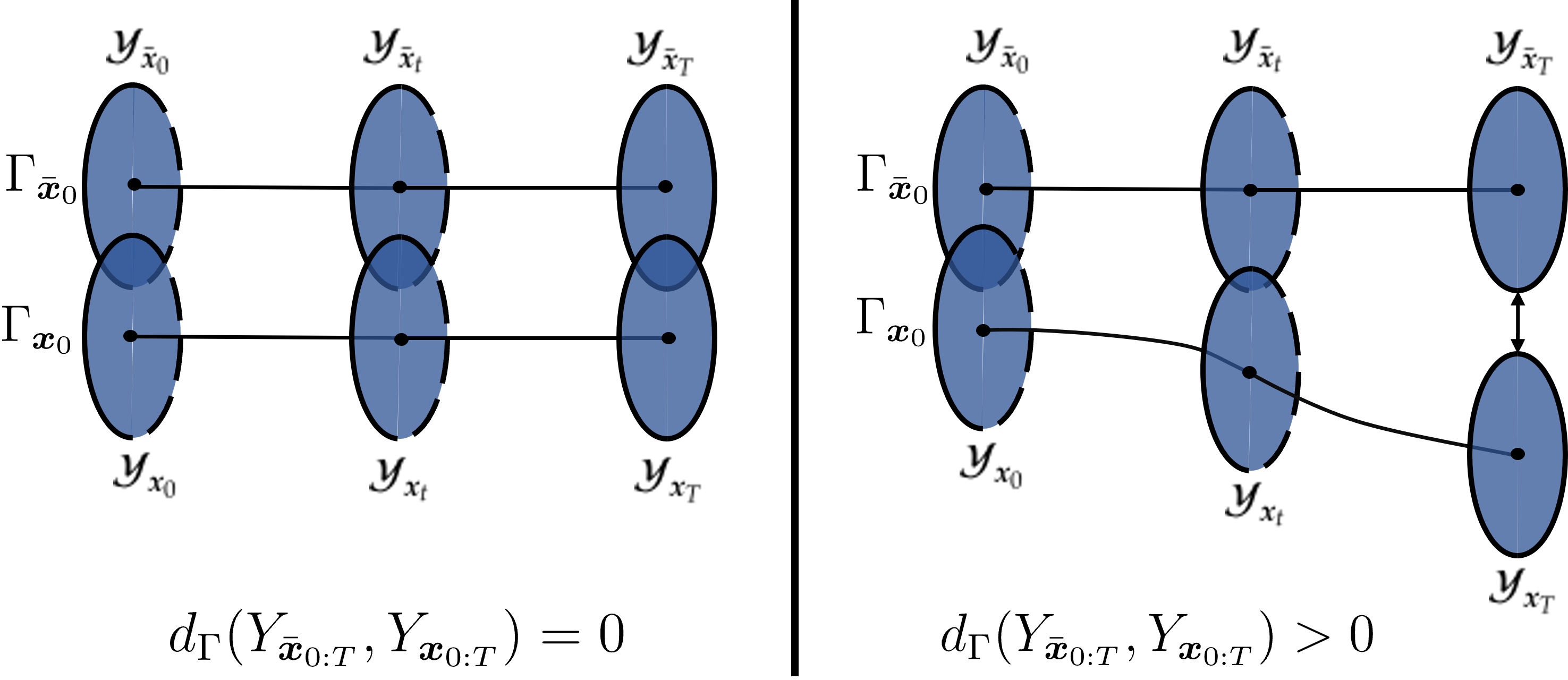}
    \caption{Representation of distinguishable and undistinguishable initial conditions using tube separation.}
    \label{fig:tubesunobs}
\end{figure}

The term ``{weakly}'' Definition~\ref{Obs} suggests that the corresponding conditions are valid for $\bar{\bm x}_0$ in a local neighborhood and not necessary 
for all $\bm x \in \mathbb{X}$. 
This definition ensures that, for finite horizon observability ($T>0$), the output tube generated by $\bm x_0$ is distinct from those generated by neighboring initial states. Fig.~\ref{fig:tubesunobs} illustrates this for $\bm x_0$, depicting unobservable  (left) and observable maps (right), respectively. In Fig.~\ref{fig:tubesunobs}, the output tubes are generated for $\bar{\bm x}_0$ and a state $\bm x$ in the neighborhood of $\bar{\bm x}_0$, with a fixed input sequence.
For the case shown on the left panel of the figure, the output tubes are indistinguishable as the separation metric is zero; as such, $\bar{\bm x}_0$ is not observable. On the contrary, the case shown on the right panel has a positive separation between the output tubes; as such, $\bar{\bm x}_0$ is observable.
An important aspect of set-valued observability of importance for path planning is the finiteness of the planned trajectory. Observability in general does not require a ``short-horizon" distinguishability condition and $T$ in Definition:\ref{Obs} can be arbitrarily large. 
A notion of short-time observability has been used for trajectory generation in deterministic output case by Alaeddini, Morgansen and Mesbahi~\cite{Alaeddini19Obs}. Here, by solving to optimize for observability in the trajectory planning problem, the intention is to search for a map where states are observable within the planning horizon. 

Therefore, an observability condition must be imposed for a fixed $T$, where any initial condition $\bm x_0 \notin B_{\delta}(\bar{\bm x}_0)$ is distinguishable from $\bar{\bm x_0}$. Here, $B_{\delta}(\bar{\bm x}_0)$ denotes a ball around $\bar{\bm x}_0$. In order to optimize for observability, the conditioning notion for observer design is defined for the set-valued analysis. A metric can then be defined to characterize observability, a term that will subsequently be embedded in an optimization objective for estimation-aware planning.  
   
\subsection*{Degree of Observability}
Now, consider an inverse map from the space of output tubes to the space of feasible state vectors, that is, $\Gamma^{-1}:Y\rightarrow \mathcal{X}$, where $Y$ represents the space of tubes and $\mathcal{X}$ the space of feasible states. If such an inverse map exists, that is, when for every tube $Y_{\bm x}$ there exists a unique initial condition $\bm x$, then conditions for finite horizon observability have been met. It is important to note here that $\Gamma^{-1}$ represents the estimation map as it takes the measurement sequence to a state vector. 
The degree of observability using a regularity metric can now be defined using this inverse map~\cite{Dontchev94Inverse,Frankowska2010,Dontchev05}.

\begin{defn}\label{def:Inv}(Inverse mapping and Metric Regularity)
Consider metric spaces $(\mathcal{X},d)$ and $(Y,d_{\Gamma})$ and the set-valued mapping $\Gamma: \mathcal{X} \rightarrow Y$. Then $\Gamma$ is \textit{metrically regular} at $(\bar{\bm x},Y_{\bar{\bm x}})\mathcal{X} \in \mathcal{X} \times Y$, if there exists constant $c>0$ and neighborhoods $U$ of $\bar{\bm x}$ and $V$ of $Y_{\bar{\bm x}}$ such that:
\begin{align}
    d(\bm x,\Gamma^{-1}(Y_{\bm x})) \leq c \, d_{\Gamma}(Y_{\bm x},\Gamma(\bm x))\quad \forall \; (\bm x,Y) \in (U \times V), \label{Eq:Gammareg}
\end{align}
where $\Gamma^{-1}(Y_{\bm x}) = \{\bm x \in \mathcal{X} | Y_{\bm x} \in Y\}$; the metrics $d$ and $d_{\Gamma}$ are defined on points and tubes, respectively, for sufficiently small neighborhoods $U$ and $V$. 
\end{defn}
The constant $c$ in Eq. (\ref{Eq:Gammareg}) is called the modulus of regularity for the set map~\cite{Frankowska2010}. When $\Gamma^{-1}$ exists, $c>0$ is a system dependent constant.

Equation (\ref{Eq:Gammareg}) ensures that a perturbation in the output tube is lower bounded by the tolerance on the corresponding initial state. For the state estimation process to be robust, one requires that the set of feasible states $\bm x$ be consistent with the output tube, and the respective distance in ${\cal X}$
scales favorably compared with perturbations in the output tube.


Lastly, given the tolerance region parameterized as $B_{\epsilon}(\bar{\bm x}_0)$, it is desired that $\bar{\bm x}_0$ be  $\epsilon$-distinguishable and  robust to perturbations in the output tubes. Therefore, it is desired that for a perturbed state $\bm x_0 \in \mathcal{A} = \{\bm x_0 | \bm x_0 \notin B_{\epsilon}(\bar{\bm x}_0)\}$, the perturbation in the output tube must be maximized. 
As such, a quantity for characterizing the worst case separation in output due to perturbation in the state is defined as the \textit{ degree of observability~\footnote{The notion of \textit{degree of observability} comes from the definition of \textit{modulus of regularity}~\cite{Dontchev05} defined for the inverse map in Definition \ref{def:Inv}. This notion of regularity is often used to analyze perturbations on set-valued maps.}}, denoted by, 
\begin{align}
D_\mathcal{O}\left(Y_{\bar{\bm x}_{0:T}}\right) &=\inf\limits_{\bm x \in \mathcal{A}} d_{\Gamma}\left(Y_{\bar{\bm x}_{0:T}},Y_{\bm x_{0:T}}\right); \label{Eq:obs_metric}
\end{align}
as such, for initial states in $\mathcal{A}$ in Eq. (\ref{Eq:obs_metric}) all trajectories are distinguishable from those initialized from $\bar {\bm x}_0$ if this quantity is positive.\footnote{For example, the set can be defined as $\mathcal{A} = \{\bm x| \lVert \bm x - \bm x_0\rVert = \epsilon\}.$} 
Here, estimation aware planning would involve maximizing $D_\mathcal{O}\left(Y_{\bm x_{0:T}}\right)$ with respect to $\Gamma$, as this would result in decreasing the sensitivity of the estimator to large perturbations in measurements. 
%
In the next section, the problem of selecting input sequences to maximize observability will thus be examined.
\section{Observability metric for discrete time output-tubes}\label{sec:Obs_thm}
As shown in Eq. (\ref{Eq:map_def}),  the nominal output map $\Gamma_{\bm x_0}$ is partially parametrized by the input sequence $\bm u_{0:T}$; the corresponding value of $D_{\mathcal{O}}$ is computed by evaluating trajectories in the region around the nominal trajectory with the same input sequence. This nominal trajectory is denoted by the state/input pair sequences $(\bar{\bm x}_{0:T},\bar{\bm u}_{0:T-1})$, for the system Eq. (\ref{Eq:sys1})-(\ref{Eq:sys2}), in the absence of output uncertainties. The perturbed trajectory is generated by sampling initial condition from the aforementioned set $\mathcal{A}$ in the previous section.

The analysis for degree of observability requires computation of distances between output tubes, which are state-dependent. In order to compute these tubes, a primary deterministic output trajectory $\Sigma_{\bar{\bm x}_0}\bm u_{0:T-1}$ is generated, and then a second trajectory $\Sigma_{{\bm x}_0}\bm u_{0:T-1}$, generated by perturbing the initial state to $\bm x_0$ such that $\bm x_0 \in \mathcal{A}$. The state-dependent primary tube and the perturbed tube are then computed by applying the tube map $\Gamma$ to ${\bar{\bm x}}_0$ and ${\bm x_0}$, respectively. These tubes are denoted by $Y_{\bar{\bm x}_{0:T}}$ and $Y_{{\bm x}_{0:T}}$. 

In order to compute the deterministic trajectories, we define the perturbed dynamics around $(\bar{\bm x}_t,\bar{\bm u}_t)$, in the absence of uncertainties, using the Taylor series expansion,
\begin{align*}
     \delta \bm x_{t+1} &= {\frac{\partial f}{\partial \bm x}}\bigg\rvert_{\bm x = \bar{\bm x}_t}
     \delta \bm x_t \quad+ \frac{\partial f}{\partial \bm u}\bigg\rvert_{\bm u = \bar{\bm u}_t}
     \delta \bm u_t \quad+ \mathcal{O}(\lVert\delta \bm x_t\rVert^2,\lVert\delta \bm u_t\rVert^2), \\
     \delta \bm y_{t} &= \frac{\partial h}{\partial \bm x}\bigg\rvert_{\bm x = \bar{\bm x}_t} 
     \delta \bm x_t ,
\end{align*}
where $\delta \bm x_t = \bm x_t -\bar{\bm x}_t,\;\delta \bm u_t = \bm u_t -\bar{\bm u}_t$ and $\delta \bm y_t = \bm y_t -\bar{\bm y}_t$ are deviations from the nominal states, inputs and outputs during the time interval $t \in [0,T-1]$. It is assumed that when the state deviation $\lVert\delta{\bm x}_t\rVert$ is small and $\lVert\delta{\bm u}_t\rVert$ is bounded, then the higher-order terms $\mathcal{O}(\lVert.\rVert^2)$ are negligible. The linearized time-varying dynamics is then,
\begin{align}
    \delta \bm x_{t+1} &= A_t\delta \bm x_t + B_t\delta \bm u_t, \label{eq:linx}\\
    \delta \bm y_{t} &= C_t\delta \bm x_t , \label{eq:liny} 
\end{align}
where $A_t = [\partial f/\partial x](\bm x_t,\bm u_t) \in \mathbb{R}^{n_x \times n_x},\; B_t = [\partial f/\partial u](\bm x_t,\bm u_t) \in \mathbb{R}^{n_x \times n_u}$ and $C_t = [\partial h/\partial x](\bm x_t,\bm u_t) \in \mathbb{R}^{n_x \times n_y}$ are the corresponding Jacobians in the Taylor series expansion. 
Here, a \textit{trust region} must be defined such that  higher-order terms are bounded and the local linearization is valid. As such, the maximum state deviation is restricted to a trust region where for a fixed control maximum deviation $\lVert\delta \bm u\rVert \leq \Delta_u$,
\begin{align}
  \lVert\delta \bm x\rVert \leq \Delta_x \implies  \mathcal{O}(\lVert\delta \bm x_t\rVert^2,\lVert\delta \bm u_t\rVert^2)\leq \lVert P^{1/2}\delta \bm x_t\rVert ;  \label{Eq:trust}  
\end{align}
where $P \succ0$ bounds maximum local state deviation and depends on the system.
Note that it also naturally follows that $\epsilon$ must be selected such that $\lVert\delta \bm x_0\lVert  = \epsilon \leq \Delta_x $.
A lower bound on the discrete-time degree of observability can be computed as defined in Eq. \ref{Eq:obs_metric}.
Here $\bm x_t$ is the trajectory generated by sampling $\bm x_0$ from $\mathcal{A}$ and propagating the dynamics with input sequence $\bm u_t = \bar{\bm u}_t \implies \delta \bm u_t = 0$. As discussed previously, the initial perturbation is set up $\epsilon$; thus, $ \lVert\delta \bm x_0\rVert = \epsilon$.
The perturbed output trajectory can be computed for the dynamics Eq. (\ref{eq:linx}-\ref{eq:liny}) as $\bm y_t = \bar{\bm y}_t + C_t(A_{t-1}\ldots A_1 A_0)\, \delta \bm x_0$. With perturbations defined on deterministic dynamics Eq. (\ref{Eq:sys1}-\ref{Eq:sys2}), the tubes can be generated with the map $\Gamma$. 
If for such $\bm x_0 \in {\cal A}$, $d_{\Gamma}\left(Y_{\bar{\bm x}_{0:T}},Y_{{\bm x}_{0:T}}\right)> 0 $, then $\bar{\bm x}_0$ is \textit{finite horizon weakly observable}. A lower bound on the degree of observability, using Eq. (\ref{Eq:setdistlb}) and Eq. (\ref{Eq:obs_metric}), can then be computed, i.e.,

\begin{align}
    D_\mathcal{O}\left(Y_{\bar{\bm x}_{0:T}}\right) 
    &= \inf\limits_{\bm x \in \mathcal{A}}\sum\limits_{t =0}^{T} d_s\left(\mathcal{Y}_{\bar{\bm x}_t},\mathcal{Y}_{\bm x_t}\right) \geq \sum\limits_{t =0}^{T} \inf\limits_{\bm x \in \mathcal{A}} d_s\left(\mathcal{Y}_{\bar{\bm x}_t},\mathcal{Y}_{\bm x_t}\right) \label{eq:aa}\\
    &= \sum\limits_{t=0}^{T}\inf\limits_{\bm x_0 \in \mathcal{A}} \; \max \; \left\{ \lVert \bar{\bm y}_t - \bm y_t\rVert - \Lambda(\mathcal{Y}_{\bar{\bm x}_t}) - \Lambda(\mathcal{Y}_{\bm x_t}),\; 0\right\} \label{eq:dobs1} \\
    &\geq \sum\limits_{t=0}^{T} \max \;\left\{\left(\inf\limits_{\bm x_0 \in \mathcal{A}}\;  \lVert C \delta \bar{\bm x}_t \rVert \right)-  \Lambda(\mathcal{Y}_{\bar{\bm x}_t}) - \Lambda(\mathcal{Y}_{\bar{\bm x}_t}) - L(\bar{\bm x}_t, \delta{\bm x}_t)\lVert \delta{\bm x}_t\rVert,\; 0\right\} \label{eq:dobs2} \\
    &\geq \sum\limits_{t=0}^{T} \max \left\{ \underbrace{s_1\epsilon\;}_{T1} \;\underbrace{\;-\; s_2 \epsilon L(\bar{\bm x}_t,s_2\epsilon)  -  2\Lambda(\mathcal{Y}_{\bar{\bm x}_t})}_{T2},\; 0\right\}:= D_\mathcal{O}^\ell(Y_{\bar{\bm x}_{0:T}}) . \label{Eq:Obs_lb2} 
\end{align}
Here, $s_1(t) = \sigma_{n}(C_tA_{t-1}\ldots A_0)$ and $s_2(t) = \sigma_{1}(A_{t-1}\ldots A_0)$ and $\sigma_{1}$ and $\sigma_{n}$ are the largest and smallest singular values of their respective matrix arguments.\footnote{Here, the lower and upper bounds on the induced norm of matrices have been used.}

The parameter $\epsilon$ in the above analysis is selected such that, for some positive $t \leq T $, the right-hand side of Eq. (\ref{Eq:Obs_lb2}) is positive.
Note that the separation term $T1$ and the uncertainty term $T2$ are independent of each other. Moreover, the separation term $T1$ depends only on the initial deviation $\epsilon$ and the dynamics of the system. It is desired to have $T1$ be bounded, as the optimization problem involves maximizing this lower bound for observability. When $A_t$ is not necessarily Schur stable, then $T1$ can grow exponentially, even for finite horizon problems, and scaling the augmented observability costs becomes an issue. Maximizing $T_1$ does improve observability as defined in Definition \ref{def:observable}, but it implies that the linearized dynamics are unstable, which is undesirable for tracking. Therefore, exploration terms must be bounded. The term $T2$ which maps how uncertainty grows over the trajectory also needs to be bounded as with unbounded uncertainty, separability of trajectories, and consequently observability, becomes infeasible; hence, Assumption (\ref{as:eigs}) is required for this analysis; the term $D_\mathcal{O}^\ell$ is now becomes well-defined. Note that for distinguishability, we need $D_\mathcal{O}^\ell>0$. One can select a large enough $\epsilon$ such that the derived lower bound is positive. Hence, for the optimization problem that will be posed shortly, the parameter $\epsilon$ acts as a tuning parameter that must be selected large enough for $D_\mathcal{O}^l$ to be positive for the initial trajectory. This is implied by stating that the system is observable for a particular map $\Gamma$ for an initial trajectory $\bar{\bm x}_{0:T},\bar{\bm u}_{0:T-1}$.

The lower bound $D_\mathcal{O}^\ell(Y_{\bar{\bm x}_{0:T}})$ on the degree of observability has been defined as a function of the state trajectory $\bm x_{0:T}$. Under the assumption on the output uncertainty, one obtains a system with set-valued outputs, defined here as \textit{observability-regular}. The system is \textit{finite horizon weakly observable} 
when $D_\mathcal{O}^\ell(Y_{\bar{\bm x}_{0:T}})>0$.

\begin{defn}\label{def:observable}
    The system Eq. (\ref{Eq:sys1})-(\ref{Eq:sys2}) with set-valued outputs defined in Eq. (\ref{Eq:ellipse}) is \textit{observability-regular} if Assumption \ref{as:eigs} holds and the system is \textit{finite horizon weakly observable}. 
\end{defn}

The main theorem for the observability measure is stated as follows:
\begin{thm}\label{thm:obs}
    Given a system with a finitely observable state $\bar{\bm x}_0$, the lower bound on the degree of observability is a concave function with respect to the state sequence $\bm x_{0:T}$ if the output map is \textit{observability-regular}.  
\end{thm}
Proof: Given $\bm x_0$ and an observability-regular uncertainty map, one can compute the lower bound on the degree of observability as given by Eq. (\ref{Eq:Obs_lb2}). Here $\Lambda$, and consequently $L(\bar{\bm x})$, are convex. Thus, Eq. (\ref{Eq:Obs_lb2}) is a summation of concave functions and, as such, concave. \\
\indent Note that this proof is structured on the condition that the uncertainty map is convex with respect to the state or has an upper-bounding envelope that is convex. Now, one can use this metric as a cost function in an optimal trajectory design problem. 
 One notes that it is in fact the lower bound Eq. (\ref{Eq:Obs_lb2}) that will be used in the guidance objective that, when maximized, improves the observability of the trajectory. Since the T1 term is constant within the trust region, maximizing the degree of observability with respect to the trajectory, only increases the tube separation by optimizing over the convex function $\Lambda$. For the optimization problem, max function in Eq. (\ref{Eq:Obs_lb2}), can be substituted out and objective for maximizing observability can be written as,
\begin{align}
     \widehat{D}_\mathcal{O}^\ell(Y_{\bar{\bm x}_{0:T}}) =  \sum\limits_{t=0}^{T} s_1\epsilon - s_2 \epsilon L(\bar{\bm x}_t,s_2 \epsilon) -  2\Lambda(\mathcal{Y}_{\bar{\bm x}_t}), \label{eq:cost1}
\end{align}
where $s_1 = (\sigma_{n}(C_tA_{t-1},\ldots A_0))$ and $s_2 = \sigma_{1}(A_{t-1},\ldots A_0)$ are constants computed over the parameters of the nominal trajectory.
Here, the convexity of $\Lambda$ and $L$ with respect to $\bar{x}_t$ ensures that when initialized with an observable trajectory (with $D_\mathcal{O}^\ell(Y_{\bar{\bm x}_{0:T}})>0$), the updated trajectory is ``more'' observable. Due to the convexity assumption, the optimization step at iteration $k$ from $\bar{\bm x}^k_{0:T}$ to ${\bar{\bm x}}^{k+1}_{0:T}$ guarantees that $\widehat{D}_\mathcal{O}^\ell(Y_{{\bm x}^{k+1}_{0:T}}) \geq \widehat{D}_\mathcal{O}^\ell(Y_{\bar{\bm x}^{k}_{0:T}})>0$. 

In the above analysis, the function $\Lambda(\mathcal{Y}_{\bm x})$ is derived from the sensor model and $L$ can be computed using Eq. (\ref{eq:lips2}). If the analytical form of $\Lambda(\mathcal{Y}_{\bm x})$ is known, then 
as stated in Assumption~\ref{as:eigs}, the convexity of these functions is a requirement for the existence of a unique trajectory; as such, the user must approximate the necessary bounds via a convex envelope function.
For a special case, when $\Lambda$ is defined by a quadratic form 
$\bm x ^TQ\bm x$ for some positive definite $Q$, 
the objective can be computed by defining a lower bound on Eq. (\ref{eq:dobs1}). 
In this case, one can write,
\begin{align}
    D_\mathcal{O}\left(Y_{\bar{\bm x}_{0:T}}\right) 
    &= \sum\limits_{t=0}^{T}\inf\limits_{\bm x_0 \in \mathcal{A}} \;  \left\{ \lVert \bar{\bm y}_t - \bm y_t\rVert - \Lambda(\mathcal{Y}_{\bar{\bm x}_t}) - \Lambda(\mathcal{Y}_{\bm x_t})\right\}\geq 0  \nonumber\\
    &\geq \sum\limits_{t=0}^{T} \;  \left\{ s_1\epsilon  - \Lambda(\mathcal{Y}_{\bar{\bm x}_t}) - \Lambda(\mathcal{Y}_{\bm x_t})\right\} \nonumber \\
    &= \sum\limits_{t=0}^{T} \;  \left\{ s_1\epsilon  - \bar{\bm x}_t^\top Q\bar{\bm x}_t - (\bar{\bm x}_t + \delta \bm x_t)^\top Q(\bar{\bm x}_t + \delta \bm x_t)\right\} \nonumber\\
    &= \sum\limits_{t=0}^{T} \;  \left\{ s_1\epsilon  - 2\bar{\bm x}_t^\top Q\bar{\bm x}_t - \sigma_1(Q)s_2\epsilon - 2\lVert\bm x_t\rVert\sqrt{s_2\epsilon})\right\} \nonumber\\
    &= \sum\limits_{t=0}^{T} \;  \left\{ (s_1\epsilon - \sigma_1(Q)s_2\epsilon)  - 2\bar{\bm x}_t^\top Q\bar{\bm x}_t  - 2\lVert\bm x_t\rVert\sqrt{s_2\epsilon})\right\} . \label{eq:cost2}
\end{align}
When the enveloping function that approximates the upper bound on the ellipsoid size is a quadratic function on $\bm x_t$, then the lower bound on the degree of observability can be approximated with Eq. (\ref{eq:cost2}) and used as the objective function.

It is also important to note, as stated earlier, the observability-maximizing trajectory could lead to a solution that is ``unstable'' for an infinite horizon problem. For the finite horizon setting, maximizing this metric will lead to a trajectory that seeks optimal state estimation, but in applications, this trajectory might conflict with task-completion constraints, such as reaching a desired final state. The observability maximization seeks to minimize uncertainty and naturally will prescribe states corresponding to low uncertainty. 
Therefore, an estimation-aware optimization problem is presented as a tool to generate a deviated trajectory with respect to some nominal trajectory that achieves the primary task. The observability metric can then be maximized under a constraint that limits the deviation from this nominal trajectory. Alternatively, one can augment the observability metric, with an existing trajectory design objective, in order to plan an estimation-aware trajectory.
These two scenarios highlight potential use cases of the proposed set-valued observability metric for trajectory design problems. Now two application cases that 
showcase these two scenarios are presented.
\section{Finite horizon Estimation-Aware Trajectory Planning}\label{sec:traj}

Offline trajectories generated for vehicular systems are often in the form of waypoints that are tracked via a lower-level feedback controller. Generally, this low-level controller has a time-constant that is smaller than timescales inherent in the planning trajectory generation.~\footnote{The open loop convergence time for low level controllers is referred to as the time constant of a control system. The time constant for the feedback control ensures that the tracking process does not incur a large lag due to control implementation.} The internal controllers used for tracking depend on state observers. It is assumed that these observers are designed without knowledge of the statistical properties of the measurement uncertainties. For the deterministic case, the state observers are guaranteed to converge to the ``true'' state when the system is observable. For the state-dependent and set-based uncertainties that considered in this paper, these state observers will generally have sub-optimal estimation performance and the system may not even be observable. Therefore, in order to optimize the tracking performance for these systems, the estimation-aware trajectory design is formulated such that the lower bound on the degree of observability is optimized.

In this section, two scenarios for generating estimation-aware trajectories are discussed.
In the first case, it is assumed that a nominal trajectory $\bar{\bm x}_{0:T}$ has been defined and user requires a improvement in estimation.
In this setting, the estimation-aware trajectory are realized via deviations to $\bar{\bm x}_{0:T}$, while satisfying the hard constraints of the original task. In the second case, the task and estimation-aware design problem are solved simultaneously as a solution to an optimization problem. In both scenarios, the main ingredient of the formulation involves maximizing Eq. (\ref{Eq:Obs_lb2}).  

\subsection{Maximizing observability for a pre-defined nominal trajectory}\label{sec:case1}

Consider the system described in Eq. (\ref{Eq:sys1})-(\ref{Eq:sys2}), where the linearization has been  given by Eq. (\ref{eq:linx})-(\ref{eq:liny}). The approximation error for this setting is negligible for the trust region defined by Eq. (\ref{Eq:trust}). Consider a nominal trajectory and control for achieving a certain task has been given in the form of $\bar{\bm x}_{0:T},\bar{\bm u}_{0:T-1}$. One then explore deviations from the nominal trajectory that maximizes observability while following strict constraints specified by the primary task. It is assumed that the enveloping function for the size of the uncertainty, $\Lambda(\mathcal{Y}_x)$, has been defined for this scenario. In practice, the state-dependent envelope is either known due to knowledge of the sensor's characteristics or estimated and stored as an oracle. Subsequently, depending on the mission constraints, an exploration parameter $\gamma>0$ must be defined that bounds the total deviation from the nominal trajectory. Hence, one can add an exploration constraint $\lVert \bm x_t - \bar{\bm x}_t\rVert_{2}\leq \gamma,\; \forall t\in[0,T]$, thereby bounding the total deviation due to exploration.~\footnote{The norms $\{1,2,\infty\}$ can also be chosen depending on the application. In this scenario, the aim is to bound the overall error with respect to the infinity norm.}

The estimation-aware planning problem for the given nominal trajectory can now be defined as \textbf{Problem 1}:

\begin{align}
    \min\limits_{ \bm u_{0:T-1}} &\quad  - D_\mathcal{O}^\ell(Y_{{\bm x}_{0:T}}) \label{Eq:obs_nom1} \\
    \text{s.t.}  &\quad \bm x_{t+1} = f(\bm x_t,\bm u_t) \nonumber \quad \mbox{for all} \; t\in[0,T-1],\nonumber \\
    &\quad \lVert \bm x_t - \bar{\bm x}_t\rVert \leq \gamma \quad \mbox{for all} \; t\in[0,T-1],\nonumber \\
    &\quad \bm x_0 = \bar{\bm x}_0, \nonumber \\
    &\quad \bm x_T = \bar{\bm x}_T. \nonumber
\end{align}
Here, the objective has been constructed as a convex function of the state sequence as defined in Eq. (\ref{eq:cost1}), but the dynamics are generally nonlinear, making the overall problem non-convex. A trust region-based sequential programming approach is now used to iteratively construct and solve a sequence of convex sub-problems. The trust region algorithm used here is a variation on the Basic Trust Region (BTR) algorithm defined in~\cite{scp}.~\footnote{For problems with non-convex inequality constraints, one can use for example, successive convexification algorithm described in~\cite{scvx}.} As such, the main problem is converted to a sequence of convex sub-problems by linearizing the dynamics and solving for optimal deviation from the nominal trajectory.
A trust region is constructed to restrict the optimization solution to be within a region around the nominal trajectory where linearization error is minimal. The trust region algorithm has been described in Algorithm~\ref{alg:BTR}, where the exploration constraint is used as a termination condition. When the exploration condition is reached for an iteration sub-problem, the trust region is reduced and the subproblem is re-instantiated. This step is repeated until trajectory updates between the iterations become negligible.
 
The trust regions are derived from bounds on linear approximations of the dynamics as defined in Eq. (\ref{Eq:trust}).For linear approximation of dynamics Eq. (\ref{eq:linx}-\ref{eq:liny}), the user defines the initial trust region bound $\Delta_{\bm x}^0$ such that bound defined in Eq. (\ref{Eq:trust}) is satisfied. Here, for the $k^{\text{th}}$ iteration, one has the linear approximation defined as $\delta \bm x^{k}_{t+1} = A_t^k\delta\bm x^k_t + B_t^k\delta\bm u^k_t,\; \forall\; t\in [0,T-1]$, for $\delta \bm x^k = \bm x^k_t - \bm x^{k-1}_t$ and $ \delta \bm u^k_t = \bm u^k_t - \bm u^{k-1}_t$. At the $k^{\text{th}}$ iteration, the linearization is performed around the sequence $\bm x_{0:T}^{k-1},\bm u_{0:T-1}^{k-1}$, and the subproblem is defined as \textbf{Problem 2} below,
\begin{align}
    \min\limits_{ \delta \bm u_{0:T-1}^k} &\quad  - D_\mathcal{O}^\ell(Y_{{\bm x}_{0:T}}^k) \label{Eq:obs_nom2} \\
    \text{s.t.}  &\quad  \delta \bm x^{k}_{t+1} = A_t^k\delta\bm x^k_t + B_t^k\delta\bm u^k_t, \nonumber \quad \mbox{for all} \; t\in[0,T-1],\nonumber \\
    &\quad  \lVert \delta \bm x^{k}_{t}\rVert_{\infty} \leq \Delta_x^k \nonumber \quad \mbox{for all} \;i\in[0,T-1],\nonumber \\
    &\quad \delta \bm x_0^k = 0, \nonumber \\
    &\quad \delta \bm x_T^k = 0. \nonumber
\end{align}
The state constraints above follow directly from \textbf{Problem 1}, ensuring that the initial conditions and final goal for the estimation-aware trajectory remain the same as that of the nominal trajectory. The exploration constraint shown in \textbf{Problem 1} Eq. (\ref{Eq:obs_nom1}) sets the termination criterion for the iterative algorithm as defined in Section~\ref{BTR}.

In order to make the problem formulation more concrete, consider the so-called
nonlinear unicycle model; in fact, the model is a Dubins car with a double integrator system rather than a conventional single integrator model~\cite{Dubin}; in order to streamline the presentation, no turning-radius constraints are enforced. Specifically, the nonlinear dynamic model is defined as,
\begin{align*}
    \dot{v}_x &= u_1 \cos\theta, \\
    \dot{v}_y &= u_1 \sin\theta, \\
    \dot{\omega} &= u_2,
\end{align*}
where $v_x,v_y$ denote the linear velocity in the $x,y$ directions, and $\theta$ is the vehicle's  orientation with respect to the world frame. The derivatives $\dot{v}_x,\dot{v}_y$ are the acceleration dynamics of the agent driven by the forward thrust $u_1$ and the
turning torque $u_2$. The output is measured as the position of the vehicle with a state-dependent uncertainty. 
The description of the vehicle dynamics for the estimation problem and linearized discrete time dynamics are discussed in Section~\ref{ap:dubin} where the system is converted into the form Eq. (\ref{eq:linx})-(\ref{eq:liny}) for the vehicle state, control input, and output given as (${\bm x},{\bm u},{\bm y}$). 
The output for the discretized system is uniformly sampled from the uncertainty set such that $\bm y \sim U(\mathcal{Y}_{\bm x})$,  where $\mathcal{Y}_{\bm x} = \{\bm y|\; C\bm x + K\lVert C\bm x - \bm y_s \rVert_2^2 + \bm r_0\}$. This candidate function has been chosen to emulate sensors that are state dependent but their uncertainty does not change arbitrarily with respect to the state. 
In more general cases, only a uniform maximum bound might properly characterize the uncertainty set. Instead, here the focus is on cases where the sensor state dependence property has been predetermined by the user and its upper bound can be approximated as a convex function. One notes that this is a valid assumption for physical systems where this model abstracts the local behavior of state-dependent sensors. For example, sensors that operate based on range-type measurements such as LiDAR often exhibit this behavior~\cite{Xu22}. The purpose of using the candidate model is to show how to generate estimation-aware trajectories. As such, the presented methodology can be realized in 
applications where the state dependency of the sensor's error is set-theoretically well defined.

\begin{figure}[t]
\centering
\subcaptionbox{Nominal and estimation-aware trajectory comparison starting at \textbf{(1,1)} ending at {\color{green}\textbf{(0,0)}} with light source at {\color{red}\textbf{(1,0)}}  \label{fig:dubin1}}
{\includegraphics[width = 0.57\linewidth]{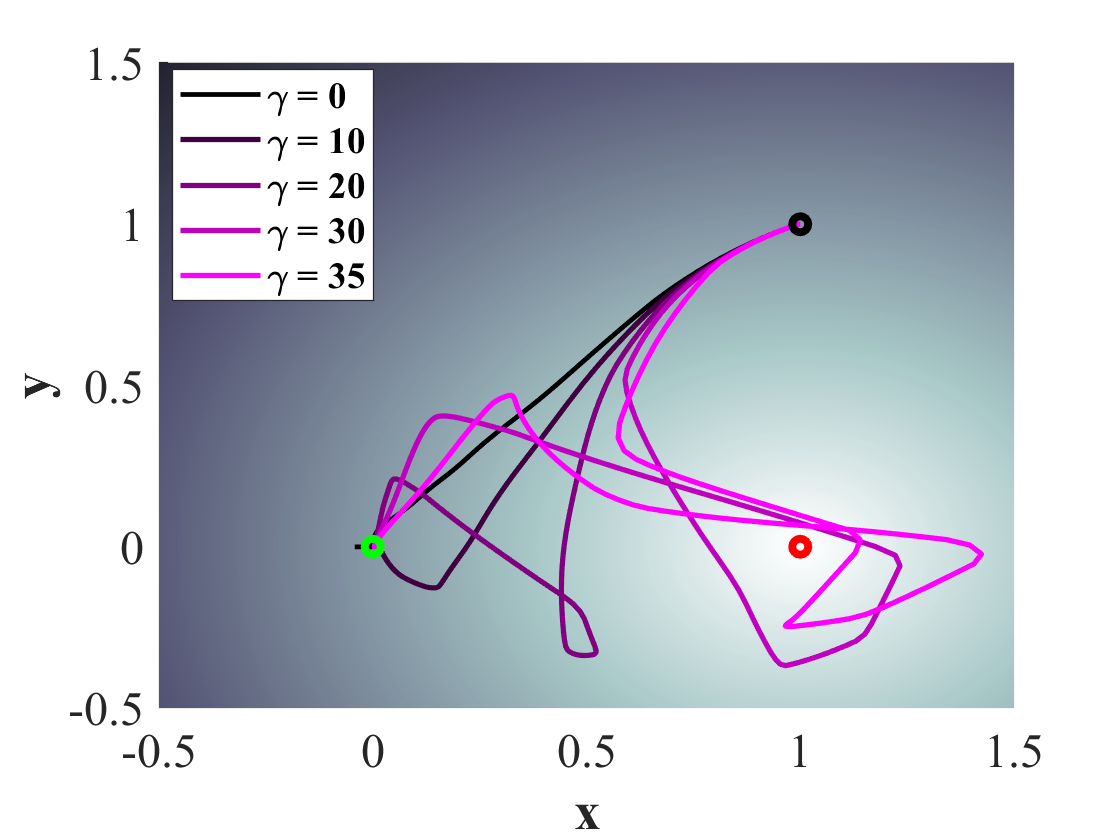}}
\subcaptionbox{Comparing the state error while tracking nominal and estimation-aware trajectories over time over 10000 roll-outs. \label{fig:dubin2}}
{\includegraphics[width=0.43\linewidth]{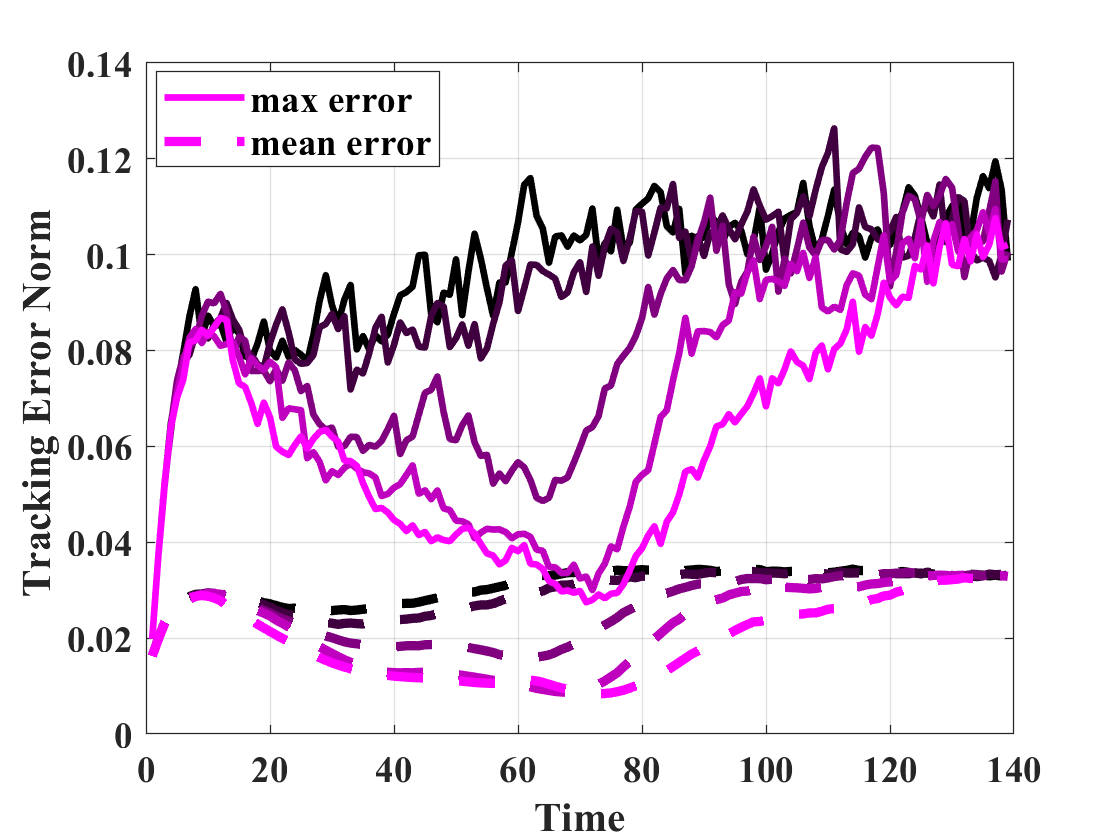}}\hfil
\subcaptionbox{Comparing the state variance across estimation-aware trajectories over time for 10000 roll-outs.\label{fig:dubin3}}
{\includegraphics[width=0.43\linewidth]{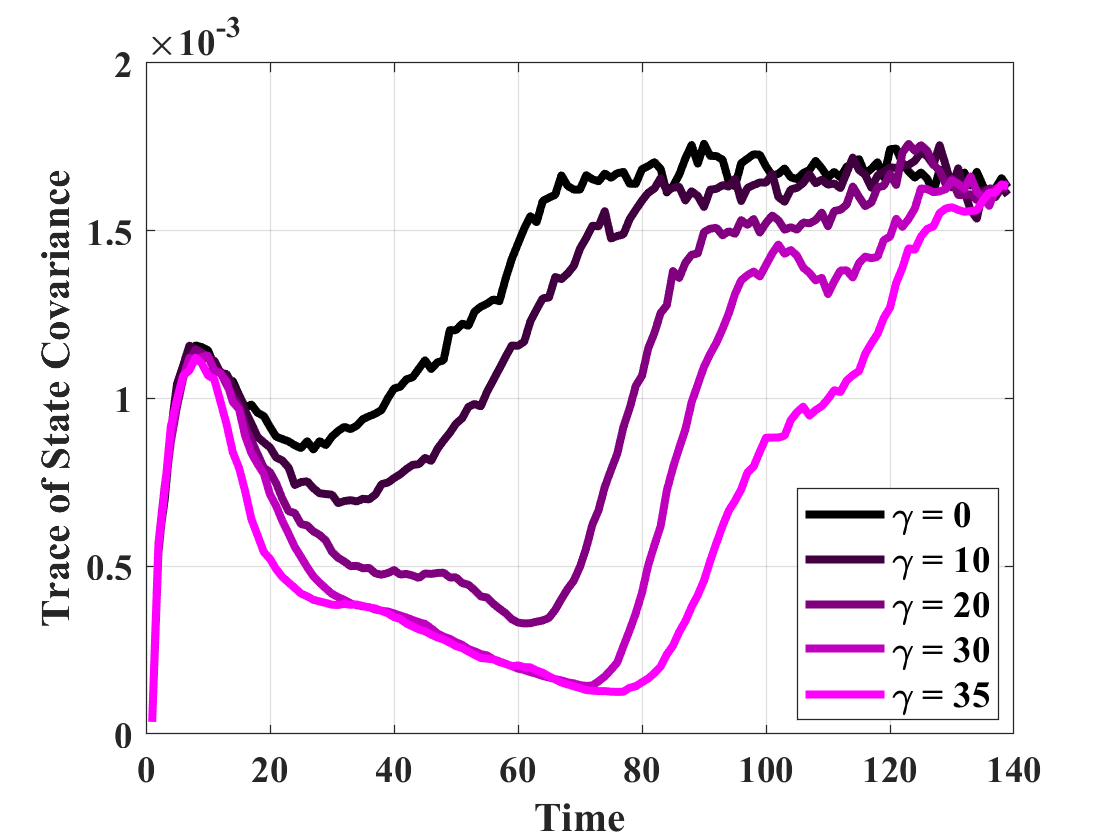}}
\caption{Comparative Analysis for Dubins car example with nominal trajectory and estimation-aware trajectories for increasing exploration bounds are presented.}\label{fig:dubin}
\end{figure}

For the uncertainty set in subsequent discussion, $\bm y_s$ denotes the location of the light source and $\bm r_0$ is the residual error that is independent of the state. This error profile is represented by the color gradient in the background of Fig.~\ref{fig:dubin1}, where the light source is shown by the red circle and the darker regions represent higher output uncertainty at the location. This measurement  model has been chosen for the problem considered not only for its practical relevance, but also as it captures a state-dependent uncertainty whose size $\Lambda(\mathcal{Y}_x)$ has a convex variation with respect to the state. The maximum size of the measurement uncertainty is characterized by a linear function of the agent's distance from the light source, making $\Lambda(x)$ convex. Since for this case Assumption~\ref{as:eigs} holds, the output map is \textit{observability-regular} due to Theorem~\ref{thm:obs}. 
The nominal trajectory for \textbf{Problem 1} is obtained by solving the primary objective where the agent/vehicle starts at position $(1,1)$ and reaches the goal at $(0,0)$. The nominal trajectory $\bm x_{0:T},\bm u_{0:T-1}$ has been characterized in order to optimize the objective $J(x,u) = (\bm x_T - \bm x_r)^\top Q(\bm x_T - \bm x_r) + \sum_{t=0}^{T-1} [(\bm x_t - \bm x_r)^\top Q(\bm x_t - \bm x_r) + \bm u_t^\top R\bm u_t]$, where $Q,R$ are positive semidefinite and positive definite matrices, respectively. 
Note that the nominal trajectory can be generated through other procedures. \textbf{Problem 2} defined in Eq. (\ref{Eq:obs_nom2}) is now solved using the BTR algorithm. 
In addition to ensuring the boundary conditions for the nominal trajectory generation,
in the estimation-aware trajectory generation, input saturation and state constraints are 
included in \textbf{Problem 2}. The exploration parameter is a design parameter for the planning problem.

The resulting estimation-aware trajectories and the estimation performance statistics are shown in Fig.~\ref{fig:dubin}.
Fig.~\ref{fig:dubin1} depicts the results of the estimation-aware optimization as a function of the exploration parameter $\gamma$ compared with the nominal trajectory. The nominal trajectory starts at location $(1,1)$ with the initial vehicle orientation of -$\pi$(facing -$x$ axis). 
The trajectory is designed such that vehicle reaches the goal $(0,0)$ at $-\pi$ orientation. This trajectory is shown in black and also corresponds to the solution of the estimation-aware planning when no exploration is allowed, i.e., when $\gamma = 0$. 
Here, the trajectories with higher exploration constants are represented with a lighter shade of pink. Fig.~\ref{fig:dubin1} shows that in order to optimize for observability, the solutions for trajectories with $\gamma>0$ "explore" regions with lower output uncertainty. The intended optimized trajectory does not prioritize optimizing for reducing output uncertainty but instead, optimizes for state uncertainty. The degree of observability metric as defined in Eq. (\ref{Eq:obs_metric}), is designed such that the size of the uncertainty for the state estimation process 
generated by an observer over this trajectory is reduced.
In order to compute the state estimation performance over the trajectories, assume that a simple lower-level state-feedback controller for the Dubins car agent exists. Since
the focus of this work is on the offline trajectory design, the design of a low-level controller for trajectory tracking is not discussed, as it is also system dependent. 
Nevertheless, the controller is designed to closely follow the waypoints with a fast convergence rate compared with the discretization time used in the system Eq. (\ref{eq:linx})-(\ref{eq:liny}).
  The agent is equipped with a Lurenberger observer $K_t$, such that a state estimate can be derived from the output measurements.~\footnote{Here a Lurenberger observer is designed for the linear approximations around the planned trajectory~\cite{Lureobs}; the corresponding poles are placed such that the observer has a fast error convergence.} The existence of this observer justifies the assumption that an inverse map $\Gamma^{-1}$ exists. Now, the state of the online agent can be represented as $\hat{\bm x}_{t+1} = (A_t - K_tC_t)\hat{\bm x}_t + B_t\bm u_t + K_t\bm y_t$, with the output uniformly drawn from the uncertainty set as $\bm y_t \sim U(\mathcal{Y}_{\bm x_t})$.~\footnote{Other estimator design approaches can be used with the assumption that the system, without state-dependent uncertainties, remains observable.} The agent now follows a planned trajectory with the stated waypoints and experiences errors due to the corresponding estimation uncertainties. The tracking error norm $\lVert \bm e_t\rVert = \lVert\bm x_t - \bar{\bm x}_t\rVert$, over the trajectory, is computed as a metric for tracking efficiency. For analysis of performance around each trajectory, a Monte Carlo simulation for $N=10,000$ trajectories (roll outs) with the same initialization and the sensor model for each case of $\gamma$ shown in Fig.~\ref{fig:dubin1}. The measurements are drawn uniformly from the set defined by $\mathcal{Y}_{\bm x}$ for each roll-out. Fig.~\ref{fig:dubin2} shows the error norm over time and the maximum error norm over time over $N$ roll-outs for each $\gamma$. The maximum error is shown by solid lines and the average error is shown by dashed lines corresponding to each estimation-aware trajectory shown with shades of pink.
Observe that the expected total error (integral of tracking error over time) is reduced with the estimation-aware trajectories and the worst-case state estimation errors are also reduced over time.   
Fig.~\ref{fig:dubin3} depicts the trace of the state covariance ($E[(\hat{\bm x} - E\hat{\bm x})(\hat{\bm x} - E\hat{\bm x})^\top]$) over $N$ roll-outs as a function of time. The trace of the variance represents sum of the deviations across all states, in this case, position, orientation and velocities. The total variance shows a significant improvement for scenarios where a higher level of exploration has been permitted. The trade-off for improving the estimation performance is an increased deviation from the nominal task-based trajectory and a possibility of reaching the boundary of the state and control constrained sets.
   
\subsection{Observability based augmented trajectory design for satellite rendezvous}\label{sec:sat}
The estimation-aware trajectory planning problem can be solved alongside the original task optimization, by augmenting the observability metric Eq. (\ref{Eq:Obs_lb2}), with a penalty term denoted as $\lambda_{\text{obs}}>0$. 
This penalty term determines how priority must be given to observability against task completion. 
In practice, the validation step must first be defined to explore the output space and determine the enveloping function given in Assumption~\ref{as:eigs}, ensuring that the observability problem is feasible. Given a convex objective function for the task optimization problem, expressed as $J'(\bm x,\bm u) = \sum_{t=0}^T c(\bm x_t,\bm u_t)$, consider its estimation-aware version via the augmented objective function as, 
\begin{align}
    \min_{\bm u_{0:T-1}} &  \quad \left(\sum\limits_{t=0}^T c(\bm x_t,\bm u_t)\right) - \lambda_{\text{obs}}D^\ell_O(Y_{\bm x_{0:T}}) \label{Eq:obs_aug}\\
    &\quad \bm x_{t+1} = f(\bm x_t,\bm u_t) \quad \forall t=0,\cdots,T-1, \nonumber  
\end{align}
where $\bm Y_{\bm x_{0:T}} = \Gamma_{\bm x_0}\bm u_{0:T-1}$, emphasizing that the observability metric implicitly depends on the input sequence. The objective function used here is defined in Eq. (\ref{eq:cost2}). The problem formulation consists of a convex objective function, nonlinear dynamic constraints, and state constraints. Is is assumed that the dynamics and constraint functions are differentiable. In order to iteratively solve this planning problem, the Successive Convexification ((SCvx)) Algorithm ~\cite{scvx} can be used. Any other programming approach that handles the nonlinear dynamics and non-convex constraints can be applied here depending on the problem.
Specifically, at each iteration $k$, the following optimal control problem is solved around the trajectory $\bm x^k_{0:T},\bm u^k_{0:T-1}$,
\begin{align}
    \min_{\bm \delta \bm u_{0:T-1}} &  \quad L^k(\delta\bm x^k,\delta\bm u^k)\label{eq:svxsat} \\
    &\quad \bm \delta \bm x_{t+1}^k = A_t\delta\bm x_t^k+B_t\delta\bm u_t^k + E\bm v_t^k \quad \mbox{for all} \quad t=0,\ldots,T-1, \nonumber  \\
    &\quad \lVert\delta \bm u_t^k\rVert_{\infty} \leq (\Delta_u)^k,
\end{align}
where $L^k(\delta\bm x^k,\delta\bm u^k) = \left(\sum\limits_{t=0}^Tc(\bm x_t^k+\delta\bm x_t^k,\bm u_t^k+\delta\bm u_t^k)\right) + \lambda_{\text{obs}}D^\ell_O + \sum\limits_{t=0}^T \lambda P(E^k_tv_t)$ and the solution to the current iteration is given by $\bm x = \bm x^k+ \delta \bm x^k,\;\bm u = \bm u^k+ \delta \bm u^k$. Here, $P$ is the exact penalty function and $\lambda$ is the penalty weight defined in~\cite{scvx}. The term $\bm v^k_t$ is an optimization variable called the \textit{virtual control buffer}. This term is used to make the problem feasible when the convexification step leads to dynamic infeasibility. Therefore,  the matrix $E$ is chosen to have full row rank in order to provide full authority for selection of $\bm v_t^k$. The penalty imposed on this term ensures that eventually, the virtual buffer converges to zero. The iterative algorithm and the trust region update policy are described in~\cite{scvx1}.

\begin{figure}[t!]
\centering
\subcaptionbox{Estimation-aware trajectories for increasing observability conditions.\label{fig:sats}}
{\includegraphics[width=0.55\linewidth]{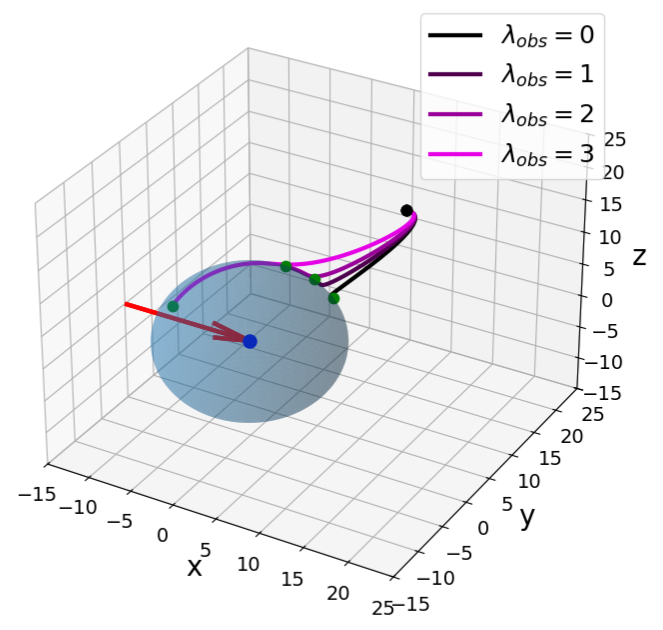}}
\subcaptionbox{Comparing tracking-error over time for for trajectories show in Fig.~\ref{fig:sats}.\label{fig:sat2}}
{\includegraphics[width=0.45\linewidth]{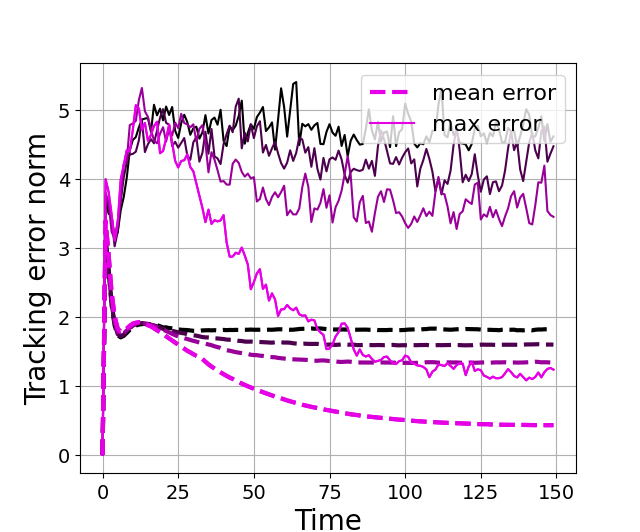}}\hfil
\subcaptionbox{Comparing the state-estimate variance for trajectories over time.\label{fig:sat3}}
{\includegraphics[width=0.45\linewidth]{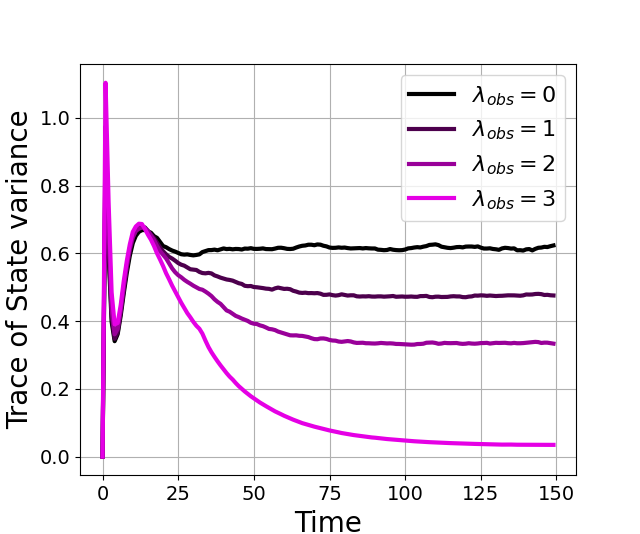}}
\caption{Comparative analysis for satellite rendezvous example. Here estimation-aware trajectories and their error analysis for increasing observability condition are presented.}\label{fig:sat}
\end{figure}

Now the trajectory generation problem for an Ego-Target rendezvous operation using the above setup can be discussed. The underlying model describes the relative dynamics between an Ego spacecraft and an uncooperative and uncontrolled Target satellite in the Hill's frame. The Target's orbit has been identified, and the Ego spacecraft is initially placed in a nearby parking orbit. The objective is for the Ego spacecraft to rendezvous with the Target while avoiding a predefined keep-out zone (shown as the blue sphere), as illustrated in Fig.~\ref{fig:sat}. 

For this scenario, its has been assumed that a camera on the Ego spacecraft captures images of the Target, and a machine learning (ML) algorithm estimates the Target’s pose (position and orientation) relative to the Ego spacecraft. The output of the ML algorithm, serves as the observed output of the system in the Target-Ego Hill's frame. A tracking controller maintains the attitude of the Ego spacecraft to ensure that the Target remains centered in the camera’s field of view. The implementation of such an ML-based attitude tracking setup, along with its integration with ML-based estimation, is discussed in~\cite{DeolePass}. The ML algorithm's architecture, implementation, and validation are detailed in~\cite{Ref:Becktor22}. 

The ML based approach has been used here as the demonstration case since it has been observed that the corresponding uncertainties are generally state-dependent. Specifically, the performance of the  ML-estimator depends on key-point identification on the Target's image. As such, the accuracy of the estimation process is directly proportional to the number of keypoints detected. The visibility of these keypoints is essential for the estimation accuracy as keypoints not illuminated by the Sun cannot be identified. As a result, the accuracy of the estimation is proportional to the illumination of Target as studied in~\cite{Ref:Becktor22} and shown in Fig.~\ref{fig:tubesunobs}. 
Moreover, for this rendezvous problem, illumination is dependent on the Ego position, as Sun's relative position is fixed during the operation. Thus, illumination is dependent on the state of the system, and therefore the output uncertainty in this case is state-dependent. This type of sensor model is also applicable to other generic image processing tools that rely on feature extraction from the image.
Note that several factors may affect the sensor accuracy, and characterizing uncertainty pertaining to every feature and parameter is non-trivial. Therefore, one needs to analyze the uncertainty as a function of the agent state while exploring the parameter space at each state. For this setup,  assume that the maximum uncertainty size can be closely approximated as a convex set with respect to the state, as described in Assumption~\ref{as:eigs}; for this example, a set defined by a quadratic function will be assumed. 
Hence, the sensor uncertainty for this case is upper bounded by the quadratic enveloping function given by $\widehat{\Lambda}(\mathcal{Y}_{\bm x}) = \bm x^TQ\bm x + p^T\bm x + r$, where $Q\in \mathbb{S}^{n_x}_+, p\in \mathbb{R}^{n_x}, r\in \mathbb{R}$. The uncertainty set $\mathcal{Y}_{\bm x}$, if not predetermined, can be approximated from the data as discussed in Section~\ref{app:sim}.

It is important to emphasize that for different scenarios, features that define state-dependency of the measurement uncertainty are distinct. Moreover, the architecture, parameters and training datasets of ML algorithms will have a significant impact on the state dependency of the output uncertainty. 
The implementation details for this scenario are discussed in Section \ref{app:sim} 


The estimation-aware trajectory design problem just described 
satisfies conditions of Theorem ~\ref{thm:obs}; as such, it can be formalized as Eq. (\ref{Eq:obs_aug}). This problem involves the synthesis of a trajectory such that the Ego spacecraft approaches the Target with some initial separation and velocity relative to the Target. The Ego spacecraft is required to approach the Target while keeping out of a zone of five meters from the Target. The design objective is formulated  such that the quadratic terms add penalty to the relative Ego approach error and fuel consumption. The observability metric augments this cost with a regularization term used as a tuning parameter. 
The overall optimization problem is thereby defined as,
\begin{align}
    \min_{\bm u_{{0:N-1}}} &  \quad \left(\sum\limits_{t=0}^T c(\bm x_t,\bm u_t)\right) - \lambda_{\text{obs}}D^\ell_O(Y_{\bm x_{0:T}}) \label{Eq:sat_obs}\\
    &\quad \bm x_{t+1} =  f(\bm x_t,\bm u_t) \quad \forall t=0,\ldots,N, &  \quad  \text{Discretized relative dynamics} \\
    &\quad ||C\bm x_t|| \geq d \quad \forall t=0,\ldots,N,&  \quad \text{Keep-out zone} \label{eq:constraint}\\
    &\quad \bm u_t \sim \mathcal{U}_{\mathcal{O}}, \quad \bm Cx_t \sim \mathcal{X}\quad \forall t=0,\ldots,N. & \quad \text{State and action constraints} \label{eq:state-action-constraint}
\end{align}
\begin{figure}[t!]
\centering
\subcaptionbox{Schematic representation of the Ego-Target Rendezvous problem.\label{fig:sat1}}
{\includegraphics[width = 0.45\linewidth]{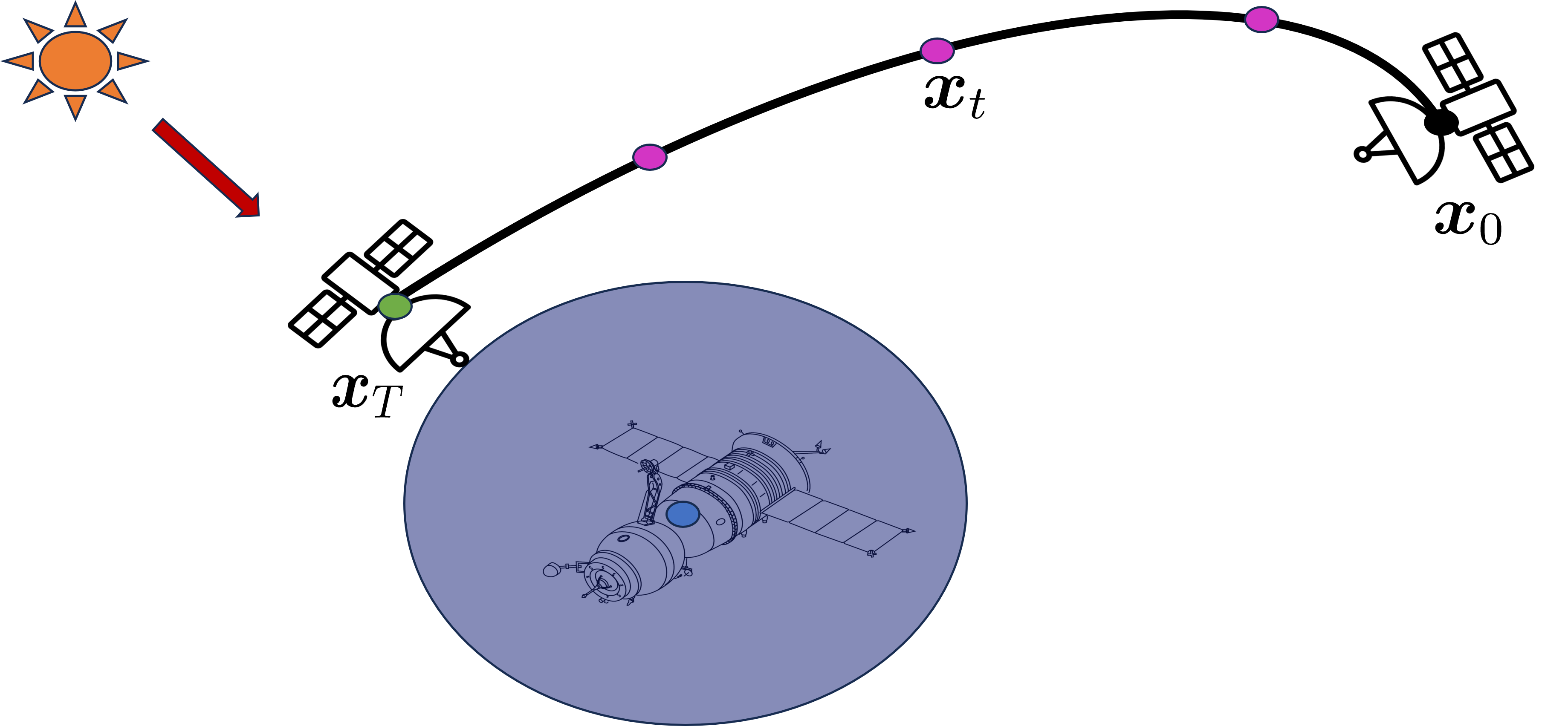}}\hfil
\subcaptionbox{Trajectory comparisons for Ego spacecraft.\label{fig:sat0}}
{\includegraphics[width=0.45\linewidth]{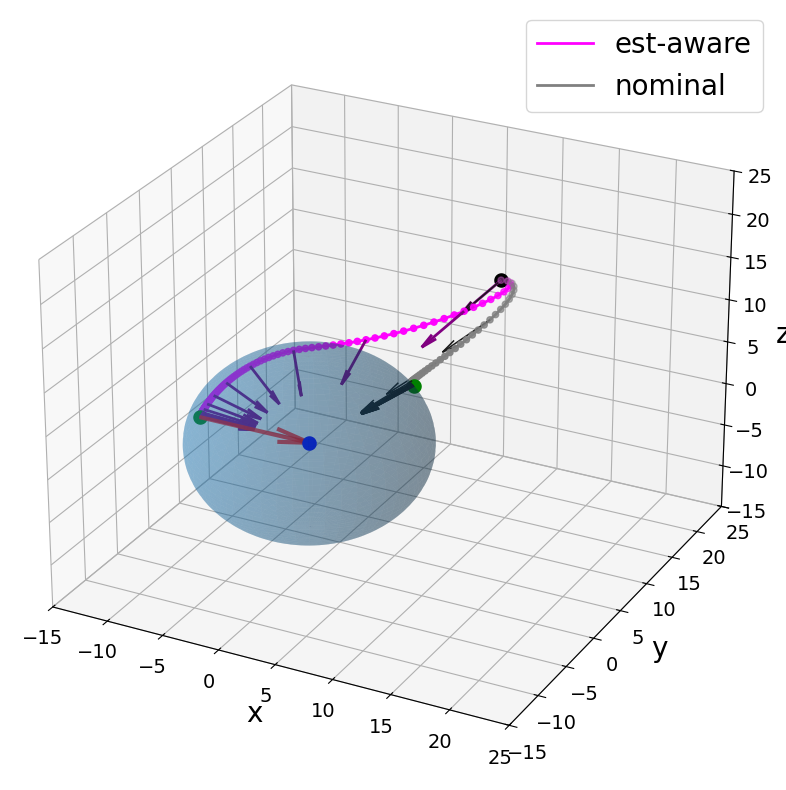}}
\caption{Ego Agent trajectory for Target rendezvous with a keep-out zone.}\label{fig:satex}
\end{figure}

 Here, the regularization term $\lambda_{\text{obs}}$ parameterizes the permissible deviation from the nominal trajectory. The objective function $J$ is designed to be convex as $l(\bm x_k,\bm u_k) = \bm x_k^\top Q\bm x_k + \bm u_k^\top R\bm u_k$ and $-D_\mathcal{O}^\ell$ are both convex. As the enveloping function is defined with a quadratic form, the cost function defined in Eq. (\ref{eq:cost2}) is used here.
 The discretized dynamics in this setting is defined in Eqs. Eq. (\ref{Eq:HCWx})-(\ref{Eq:HCWz}). The keep-out zone constraints shown here are non-convex, and as such, are linearized at each iteration of the SCvx approach. 
 Using the iterative SCVx approach, the convex subproblems as described as Eq. (\ref{eq:svxsat}). The implementation details are described in Section~\ref{sec:scvx_imp}.

Representative trajectories for this planning problem formulated as Eq. (\ref{Eq:sat_obs})-(\ref{eq:state-action-constraint}), are shown in Fig.~\ref{fig:sat0}. The proposed solution strategy thereby obtain an estimation-aware trajectory such that the Ego camera points to the Target with the optimum illumination angle.
Fig.~\ref{fig:sat1} depicts the schematic for the problem where the Target is shown in the center within the keep-out zone in blue and the Ego agent following the proposed trajectory. The Sun direction is shown by the red arrow. Note that aligning with the Sun-angle reduces the instantaneous output uncertainty for the keypoint-based ML estimator. 
Fig.~\ref{fig:sat0} depicts the derived solutions for the nominal trajectory which is designed without the augmented observability cost compared with the estimation-aware trajectory that is the solution to Eq. (\ref{Eq:sat_obs})-(\ref{eq:state-action-constraint}). The Ego trajectories for the nominal case is shown in gray, where the Target relative position is the origin and its keep-out zone is shown by the blue sphere. This is the best trajectory for a rendezvous scenario. The Ego trajectory is initialized at the location $[7,20,10]^\top$m shown by the black dot in Fig~\ref{fig:sat0}, and an initial velocity of $[10,-5,-5]^\top$m/s relative to the Target shown in green at the origin. The Ego reaches the goal location shown in green with relative position vectors shown by the arrows on the trajectory. The estimation-aware trajectory is shown in pink, where the final Ego position is such that the relative Ego position achieves the best viewing angle during the course of its trajectory.

The choice of the regularization term $\lambda_{\text{box}}$ signifies the weighting of the planned trajectory between estimation and rendezvous tasks. The estimation performance over the nominal and estimation-aware trajectories and the progression of uncertainty for different cases of $\lambda_{\text{box}}$ are shown in Fig.~\ref{fig:sat}. For these simulations, a simple low-level controller that tracks the waypoints on the designed trajectories have been chosen.
The trajectories are followed by tracking the waypoints $\bm x_{0:T}$ generated for each $\lambda_{\text{obs}}$. As also mentioned previously, a Lurenberger state observer has been designed for the discrete time system and tracking error norm is computed at each waypoint. Similar to the previous scenario, the online system is $\hat{\bm x}_{t+1} = (A-LC)\hat{\bm x}_{t+1} + B\bm u_t + L\bm y_t$, where the output is drawn uniformly from the uncertainty set given by $\bm y_t \sim U(\mathcal{Y}_{\bm x_t})$. The tracking error norm can be computed by $\lVert\bm e_{t}\rVert = \lVert\hat{\bm x}_t - \bm x_{t}\rVert$. 
Now, for analyzing the tracking performance, a Monte Carlo simulation with $N=1000$ roll-outs for the system for trajectory designed by each case $\lambda_{\text{obs}}=[0,1,2,3]$ has been developed. Fig.~\ref{fig:sat2} depicts the error norm over time $\lVert e_t\rVert$ for each trajectory. 
The dashed lines show the mean error for all roll-outs for each trajectory and the solid lines show the maximum error norm over each roll-out. Here we can observe that the full state error is significantly reduced for the estimation-aware trajectories and it improves as more importance is given to the estimation performance. Note that this is not just the measurement error but the state estimation error, which highlights that tracking for the velocity state improves as well (not simply the static position outputs). Next, the trace of the state covariance $\mbox{Trace}\left(E[(\bm x_t - E[\bm x_t])(\bm x_t - E[\bm x_t])^\top]\right)$, at each time for $N=1000$ roll-outs, has been  computed for each $\lambda_{\text{obs}}$. In Fig.~\ref{fig:sat3}, the progression of the trace of the state covariance over time has been shown. Here, one observes  that the total state variance is reduced for the estimation-aware trajectories as more importance is given to the measure of observability. Moreover, since the Rendezvous problem is a free-final state problem, we observe that the Ego spacecraft reaches different final states based on $\lambda_{\text{obs}}$.

For the two scenarios presented in this section, estimation-aware trajectories have been designed using the proposed degree of observability metric. As seen in Figs.~\ref{fig:dubin} and \ref{fig:sat}, the state estimation process significantly improve using the prescribed methodology. 
This proposed approach can be interpreted as improving the sensor conditioning, where the uncertainty in the measurement is mapped to a smaller uncertainty in the state estimation process. Moreover, along the estimation-aware trajectory, the lower-level state feedback controller would generally require less energy due to smaller feedback errors, thus improving the tracking performance. 
The proposed approach also has application in multi-agent settings, where observations of other agents' behavior for the purpose of coordination has a set-valued uncertainty.
%
\section{Conclusions}

In this paper, an observability-based metric for set-valued measurements induced by state-dependent errors was presented. This metric can be utilized to design estimation-aware trajectories for set-valued uncertainty measurements when a certain convexity assumption holds. Specifically, the output uncertainty map must be convex in its size with respect to the state Eq. (\ref{Eq:maxellipse}). In practice, when this map is not smooth, an enveloping convex function must be approximated through a validation step to establish a tight upper bound on the uncertainty. Within this setup, the manuscript demonstrates that when the output maps are \textit{observability-regular} (Definition \ref{def:observable}), a lower bound on the \textit{degree of observability} Eq. (\ref{Eq:obs_metric}), can be defined. Under reasonable assumptions on the state dependency of the uncertainty set, this lower bound is a concave function. The paper then proceeded to show that an optimization problem can be formulated to maximize the degree of observability, thereby generating estimation-aware trajectories.

This work has also illustrated the utility of the proposed observability metric in formulating an augmented cost function for trajectory design problems. A case study was considered, in which an ML-based state estimator, implemented by an Ego spacecraft, tracks an uncooperative Target satellite. The proposed procedure constrains the estimator as a state-dependent set-valued map. It has been shown that the proposed metric can generate a tracking trajectory for the Ego spacecraft, resulting in improved estimation accuracy. Furthermore, the paper has demonstrated that, given a nominal trajectory, a neighboring estimation-aware trajectory can be designed by maximizing the observability metric lower bound while imposing exploration constraints in the neighborhood of the nominal trajectory. This case was presented with an example using Dubins car model.

This proposed framework has been formulated for a general class of differentiable nonlinear systems with a set-valued output uncertainty. The optimization problems that use the observability metrics are defined for locally linearizable systems where the size of the output uncertainty is globally convex. The proposed method is effective for designing offline estimation-aware trajectories for observability-regular maps defined in the neighborhood of a nominal trajectory.

Further investigations into output uncertainties with non-ellipsoidal bounds would be a natural extension of this work. The framework can also be extended to cases where the uncertainty size can only be locally approximated as a convex function with respect to the trajectory state, while the optimization problem is solved sequentially. Additionally, further exploration of connections to other uncertainty classes in state-dependent settings and the use of a metric that quantifies output perturbation could follow as a comparative study.

Another facet of this problem is its interpretation as a dynamic sensor placement problem, where sensor trajectories can be designed to enhance an agent's state estimation. If a sensor map is provided and the sensor can be dynamically positioned around the agent it is tracking, the observability metrics can be employed to design an estimation-aware trajectory for the sensor-agent pair. Moreover, a multi-agent sensor placement scenario can be constructed, where sensors interact to maximize combined agent information as they are dynamically positioned to optimize the agent's observability.

Finally, for real-time implementation, a model predictive approach can be adopted, utilizing an updated parametrization of the output uncertainty to optimize the trajectory for estimation. This approach is applicable when online estimation begins to degrade, and the uncertainty map and model information are known.
The code for examples discussed in this paper is available  \href{https://github.com/Rainlabuw/Obs_aware_opt/tree/main}{here}.

\section*{Acknowledgments}
The authors thank their various discussions 
on estimation-aware trajectory design with Spencer Kraisler, Shahriar Talebi, Beniamino Pozzan, Saptarshi Bandyopadhyay, Amir Rahmani, and Nick Andrews. The research of the authors has been supported by NSF grant ECCS-2149470 and AFOSR grant FA9550-20-1-0053.

\section{Appendix}

Here we gather the complementary background pertaining to the main results presented in this work, and discuss the relevant implementation details.

\subsection{System dynamics for Dubins Car}\label{ap:dubin}
The nonlinear full state Dubins car dynamics is defined as,

\begin{align*}
    \dot{p_x} &= v_x,\\
    \dot{p_y} &= v_y,\\
    \dot{\theta} &= \omega,\\
    \dot{v_x} &= u_1 \cos\theta, \\
    \dot{v_y} &= u_1 \sin\theta, \\
    \dot{\omega} &= u_2, \\
\end{align*}
where the agent's state is defined as $\bm x = [p_x,p_y,\theta,\dot{p_x},\dot{p_y},\dot{\theta}]^\top$ and inputs for the forward thrust and turning torque are designated as $\bm u =[u_1,u_2]^\top $, where $p_x,p_y$ are the agent's position coordinates in 2D and $\theta$ is the heading angle in the world frame. The corresponding velocities are given by $v_x$ and $v_y$ and the angular rate by $\omega$. The output is measured by a position $\bm y = [p_x,p_y]^\top$ sensor that has an uncertainty dependent on the light intensity. 
The discretized dynamics at step $dt$ are represented as Eq. (\ref{eq:linx})-(\ref{eq:liny}), where the system matrices are,
\begin{align*}
    A(\bm x,\bm u) &= \left[\begin{array}{cccccc}
         1 & 0 & -\frac{dt^2u_1\sin\theta}{2} & dt & 0 & -\frac{dt^3u_1\sin\theta}{6}   \\
        0 & 1 & \frac{dt^2u_1\cos\theta}{2} & 0 & dt & \frac{dt^3u_1\cos\theta}{6}   \\
        0 & 0 & 1 & 0 & 0 & dt   \\
        0 & 0 & -dt u_1\sin\theta & 1 & 0 & -\frac{dt^2u_1\sin\theta}{2} \\
        0 & 0 & dt u_1\cos\theta & 0 & 1 & \frac{dt^2u_1\cos\theta}{2} \\
        0 & 0 & 0 & 0 & 0 & 1 \\
    \end{array}\right],\;
    B(\bm x,\bm u) = \left[\begin{array}{cc}
         \frac{dt^2\cos\theta}{2} & -\frac{dt^4u_1\sin\theta}{24}\\
        \frac{dt^2\sin\theta}{2} & \frac{dt^4u_1\cos\theta}{24}\\
        0 & \frac{dt^2}{6}\\
        dt \cos\theta &  -\frac{dt^3u_1\sin\theta}{6} \\
        dt \sin\theta &  \frac{dt^3u_1\cos\theta}{6} \\
        0 & dt \\
    \end{array}\right],\\
    C(\bm x,\bm u) &= \left[\begin{array}{cccccc}
        1 & 0 & 0 & 0 & 0 & 0   \\
        0 & 1 & 0 & 0 & 0 & 0   \\
    \end{array}\right].
\end{align*}

\subsection{Trust Region Algorithm Implementation}{\label{BTR}}
\begin{algorithm}
\caption{Trust Region algorithm}\label{alg:BTR}
\begin{algorithmic}
\State Initialize: Nominal trajectory ($\bar{\bm x}^0_{0:T},\bar{\bm u}^0_{0:T-1})$, Exploration Parameter:$\gamma$, Trust region : $ \Delta_x^0$.
\State Set $k=0$, $0\leq l\leq 1$, $r=0, \alpha$ 
\While{$\lVert{\bm x}^{k+1} - { \bm x}^k\rVert_2 \geq \varepsilon$} \Comment{Termination Condition}
\While{$\lVert \bar{\bm x}^0 - {\bm x}^k\rVert_2 \leq \gamma$}
\State At iteration $k$ Solve :  \textbf{Problem 2}Eq. (\ref{Eq:obs_nom2}) for ($\delta \bm x^k,\delta \bm u^k$) for $\Delta_x = \alpha*l^r$ 
    \State $\bm x_{0:T}^k \gets \bm x_{0:T}^{k-1} + \delta \bm x^k$
    \State $k \gets k+1$
\EndWhile
\If{$\lVert \bar{\bm x}^0 - {\bm x}^k\rVert_2 \geq \gamma$} \Comment{Exploration violated, update trust region}
\State discard update $\delta \bm x^k$
\State $k \gets k$
\State $r \gets r+1$
\EndIf
\EndWhile
\end{algorithmic}
\end{algorithm}

\begin{figure}[t]
\centering
\subcaptionbox{Total deviation from nominal condition for each exploration conditions.  \label{fig:btr1}}
{\includegraphics[width = 0.48\linewidth]{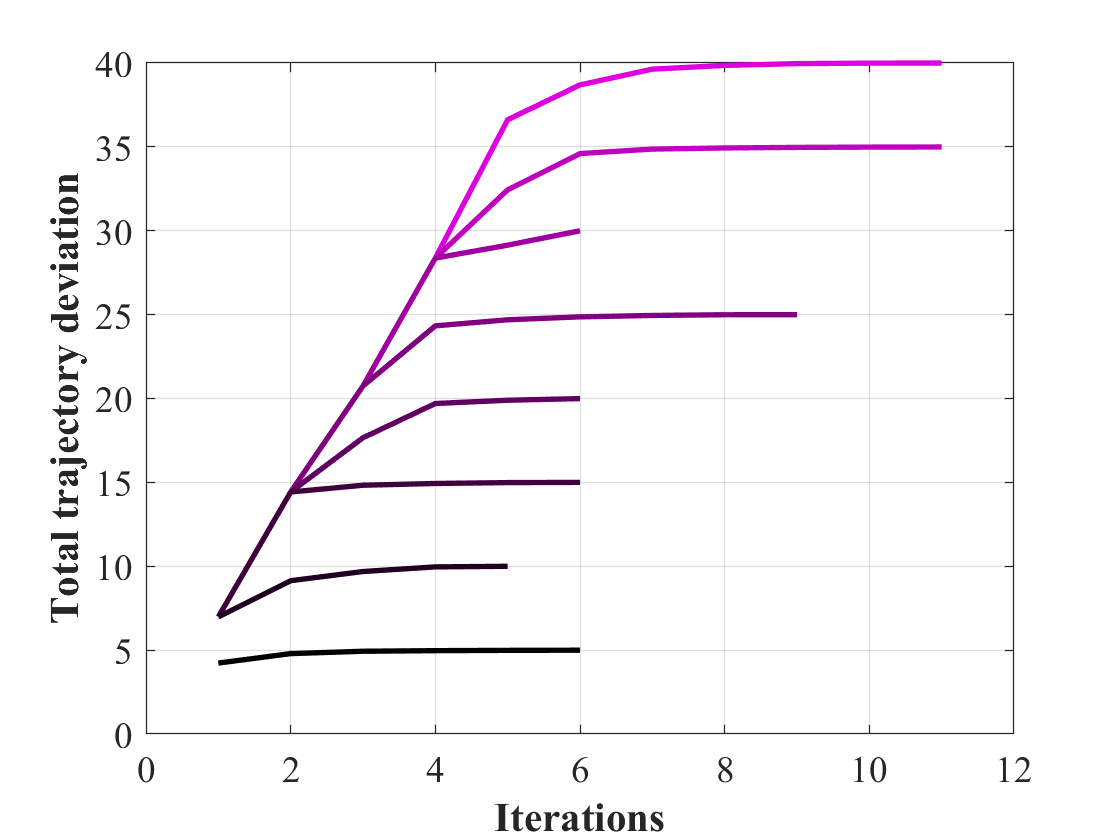}}
\subcaptionbox{Objective function growth per iteration for each exploration condition. \label{fig:btr2}}
{\includegraphics[width=0.48\linewidth]{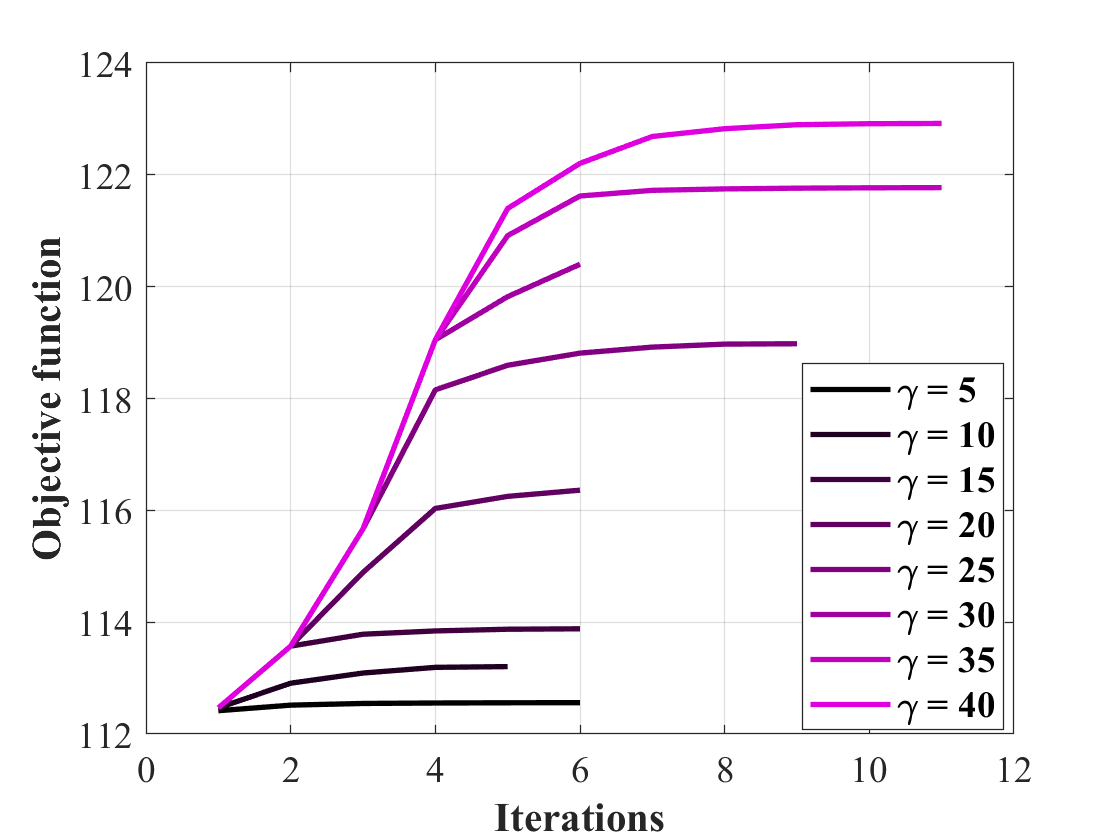}}
\caption{Trust region algorithm analysis for different observability exploration conditions used as termination conditions.}\label{fig:BTR}
\end{figure}
The algorithm we have implemented in this paper for the trajectory generation problem is based on a sequential method that involves a trust region, for which the convex subproblem is feasible. The setup is described in \textbf{Algorithm}~\ref{alg:BTR}. For the dynamics defined in the previous section, we set a trust region by setting $\lVert\delta \bm x\rVert_\infty \leq \alpha*l^r$; $l$ denotes the learning rate set to the value of $0.5$. The parameter $\alpha$ must be set such that for $r = 0$, the initial trust region $\alpha*l^r\leq \Delta_x^0$ satisfies Eq. (\ref{Eq:trust}). The trust regions shrink as iterations increase, therefore dynamic constraints are always satisfied. The trajectory at the first iteration is initialized with the nominal trajectory; at each subsequent iteration \textbf{Problem 2} Eq. (\ref{Eq:obs_nom2}) is solved in MATLAB using CVX~\cite{cvx}. We present the analysis for the performance of Algorithm~\ref{alg:BTR} in Fig.~\ref{fig:BTR}. The exploration condition $\lVert \bar{\bm x} - \bm x\rVert_2\leq\gamma$ represents the total deviation from the nominal trajectory. We observe that the trust regions updated by the algorithm ensure that the observability maximization does not violate the exploration constraints as seen in Fig.~\ref{fig:btr1}. For case of estimation aware trajectory shown in Fig.~\ref{fig:dubin}, we present the corresponding $\gamma$ evolution and degree of observability progression as a function of algorithm iterations. The trust region algorithm reduces maximum deviation as exploration constraints are reached. The shades of pink s how increasing exploration $\gamma$ from nominal trajectory on x-axis in Fig.~\ref{fig:btr1}. The $\gamma$ values for each case of estimation-aware trajectory are shown on the right in Fig.~\ref{fig:btr2}.    
The algorithm is terminated when the maximum exploration condition is reached and the trust region has been reduced below a threshold where no more deviations are feasible. Fig.~\ref{fig:btr2} shows that the degree of observability is increased by allowing more exploration. As seen in Fig.~\ref{fig:dubin}, a higher exploration limit leads to an improved estimation process. This comes at a cost as the system potentially reaches the boundary of the set for the control and state constraints; depending on the mission, this may not be desirable. Therefore, it is judicious to increase the observability metrics while monitoring their potential undesirable implications on the synthesized trajectory.

\subsection{System Dynamics for Satellite Rendezvous}\label{Sat}

The relative (translational) dynamics for the Ego spacecraft are specified in the Hill's frame of the Target. Specifically, the Clohessy–Wiltshire equations~\cite{Hills} are used to define the dynamics for the position ($x,y,z$) (coordinates in the Hill's frame) and velocity ($\dot{x},\dot{y},\dot{z}$) as,
\begin{align}
    \ddot{x} &= 3n^2x + 2n\dot{y} + u_x, \label{Eq:HCWx}\\
    \ddot{y} &= -2n\dot{x} + u_y, \label{Eq:HCWy}\\
    \ddot{z} &= -n^2z  + u_z, \label{Eq:HCWz}
\end{align}
where $n = 0.00113 \mbox{s}^{-1}$ denotes the mean motion of Target spacecraft circular orbit and $u$'s are the control inputs in each corresponding axis.
Consider now the position and velocity vectors $\bm p = [x,y,z]^\top \in \mathbb{R}^3$ and $\bm v = [\dot{x},\dot{y},\dot{z}]^\top \in \mathbb{R}^3$. First, redefine the state vector as $\bm x = [\bm p^\top,\bm v^\top]^\top \in 
\mathbb{R}^6$ and inputs as $\bm u = [u_x,u_y,u_z]^\top \in 
\mathbb{R}^3$. The output of the system measures the relative state of the Ego with respect to the Target; hence, the nominal output is defined as $\bm y = \bm p \in \mathbb{R}^3$. The problem defined in Eq. (\ref{Eq:sat_obs}) can now be solved using the discretized linear dynamics of the form Eq. (\ref{eq:linx})-(\ref{eq:liny}) for the system defined in Eq. (\ref{Eq:HCWx})-(\ref{Eq:HCWz}). The system matrices for the model $\bm x_{t+1} = A\bm x+ B\bm u;\; \bm y_t = C\bm x_t$ are defined as,
\begin{align*}
    A = \left[\begin{array}{cccccc}
         0 & 0 & 0 & 1 & 0 & 0   \\
         0 & 0 & 0 & 0 & 1 & 0   \\
         0 & 0 & 0 & 0 & 0 & 1   \\
         3n^2 & 0 & 0 & 0 & 2n & 0   \\
         0 & 0 & 0 & -2n & 0 & 0   \\
         0 & 0 & -n^2 & 0 & 0 & 0   \\
    \end{array}\right],\quad B = \left[\begin{array}{ccc}
         0 & 0 & 0    \\
         0 & 0 & 0    \\
         0 & 0 & 0    \\
         1 & 0 & 0    \\
         0 & 1 & 0    \\
         0 & 0 & 1   \\
    \end{array}\right],\quad C = \left[\begin{array}{cccccc}
         1 & 0 & 0 & 0 & 0 & 0    \\
         0 & 1 & 0 & 0 & 0 & 0   \\
         0 & 0 & 1 & 0 & 0 & 0   \\
    \end{array}\right].
\end{align*}

\subsection{SCvx implementation details}\label{sec:scvx_imp}

The SCvx algorithm converts the problem Eq. (\ref{Eq:sat_obs}) into a sequence of convex subproblems of the form,
\begin{align}
    \min_{\bm \delta \bm u_{0:T-1}} &  \quad L^k(\delta\bm x^k,\delta\bm u^k)\label{eq:svxsat2} \\
    &\quad \bm \delta \bm x_{t+1}^k = A_t\delta\bm x_t^k+B_t\delta\bm u_t^k + E\bm v_t^k \quad \forall t=0,\cdots,T-1, \nonumber  \\
    &\quad s(x^k_t) + S^k\delta \bm x_t^k -s^k_t \leq 0, \\
    &\quad \lVert\delta \bm u_t^k\rVert_{\infty} \leq (\Delta_u)^k.
\end{align}

\begin{figure}[t!]
\centering
\subcaptionbox{Trust Region size updates per iteration.\label{fig:s1}}
{\includegraphics[width = 0.45\linewidth]{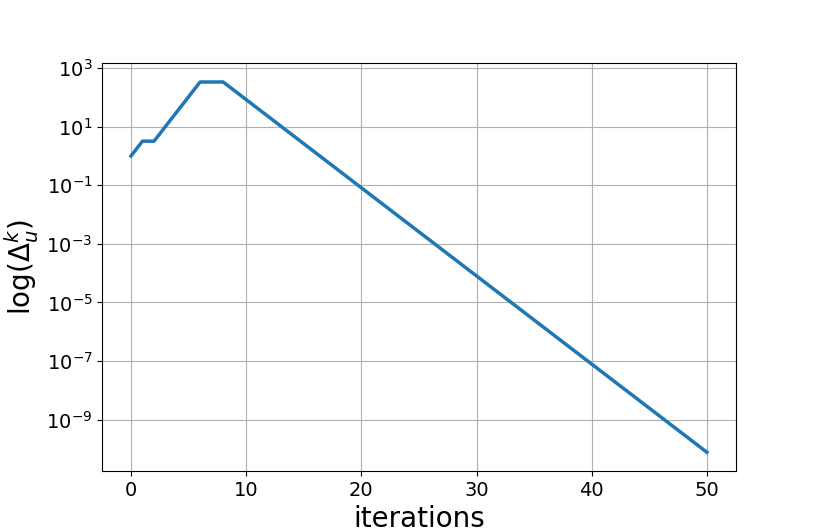}}
\subcaptionbox{Objective function update per iteration.\label{fig:s2}}
{\includegraphics[width=0.45\linewidth]{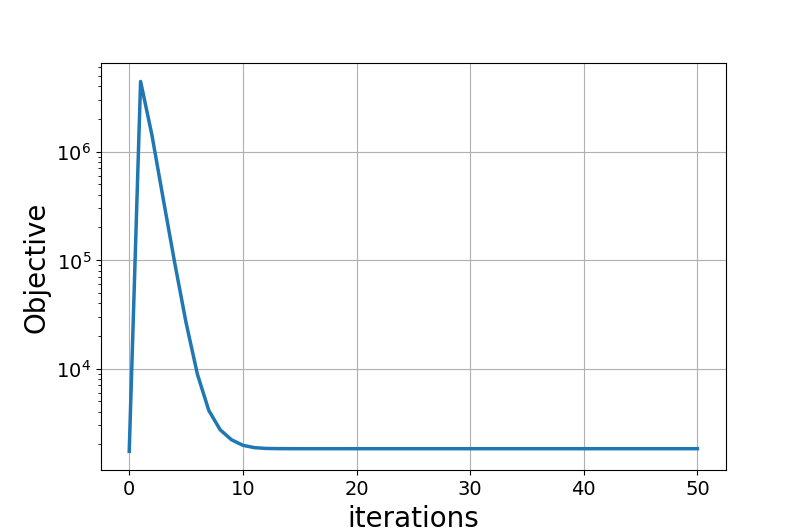}}
\subcaptionbox{Virtual control buffer evolution.\label{fig:s3}}
{\includegraphics[width = 0.45\linewidth]{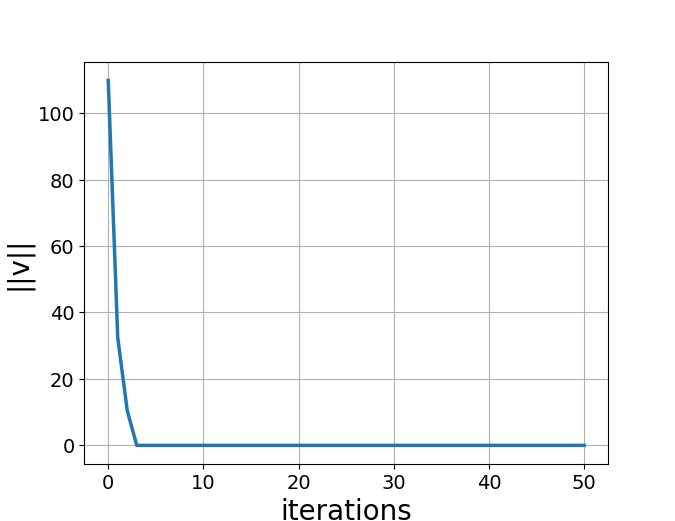}}
\subcaptionbox{Constraint violation term evolution.\label{fig:s4}}
{\includegraphics[width=0.45\linewidth]{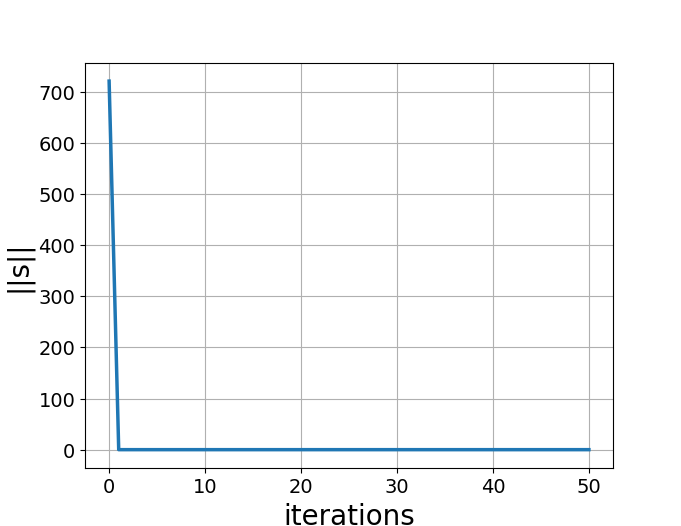}}
\caption{SCVx implementation statistics for estimation-aware Ego-Target Rendezvous problem.}\label{fig:satstat}
\end{figure}

Here $L^k(\delta\bm x^k,\delta\bm u^k) = \lambda_1\left(\sum\limits_{t=0}^Tc(\bm x_t^k+\delta\bm x_t^k,\bm u_t^k+\delta\bm u_t^k)\right) + \lambda_2\lambda_{\text{obs}}D^\ell_O + \lambda_3\sum\limits_{t=0}^T P(E^k_tv^k_t,s^k_t)$, where $\lambda$'s are the tuning penalty weights. The keep-out zone constraint in Eq. (\ref{eq:constraint}) for distance $d$ is written as $s(\bm x) = d^2 - \lVert C\bm x \rVert^2$ with the corresponding Jacobian is $S  = \partial s/\partial \bm x = -2C^\top C\bm x$. The penalty functions on the virtual buffer $v^k_t$ and the constraint violation $s^k_t$ are in terms of the $1$-norm. Here, the SCvx algorithm is applied as described in~\cite{scvx1}.
The penalty terms are chosen as $\lambda_1 = 10^{-1}$, $\lambda_2 = 10^{-1}$ and $\lambda_3 = 10^{3}$. The implementation results for the scenario shown in Fig.~\ref{fig:sat1} and the progression of trust regions and the objective function are depicted in Fig.~\ref{fig:satstat}. 
The trust region is updated to accommodate for a higher trajectory deviation and then converges as shown in Fig.~\ref{fig:s1}. The evolution of the objective function can be seen in Fig.~\ref{fig:s2}. The objective function grows initially as the penalty terms for the state constraints are high due to the initialized trajectory. Over subsequent iterations, the objective function converges. The heavy penalization of constraints violation terms ensures fast convergence of buffer terms $v$ and $s$, as seen in Figs.~\ref{fig:s3} and \ref{fig:s4}.   


\subsection{Simulation setup}\label{app:sim}

The validation experiments and online analysis are performed in a visual simulator. A closed-loop end-to-end pipeline has been developed to simulate high-fidelity spacecraft dynamics and produces photorealistic high-resolution images for the ML-based state estimator. The rendering engine used here is built as a ``level'' in Unreal Engine software, where a Low-Earth orbital environment with realistic lighting and textures are simulated.
The detailed architecture of the CNN and integrated simulation pipeline has been described in the previous work~\cite{Ref:Becktor22}.
The filtered state estimates generated by applying the perception map to the images are then used as feedback for the tracking algorithms. The simulator is a modular setup where we can select appropriate astrodynamics simulators or generate the dynamics in Python. The offline trajectory planning, which is the main optimization problem examined in this paper, has been developed using CVXPY in Python~\cite{diamond2016cvxpy}.

\subsection{Validation experiment Design}\label{sec:experiment}
In order to provide robustness guarantees, uniform error bounds must be provided for the perception map. A method similar to the procedure suggested by Dean {\em et al.}~\cite{dean2019robust} can be adapted, where they generate a slope-dependent model of the estimate error near any training point. However, in this work, the uncertainty is modeled as a state-dependent ellipsoidal set. As such, with the assumption that the output uncertainty with respect to each state is bounded, bounds are estimated from the data. 
\begin{figure}
    \centering
    \includegraphics[width=0.75\linewidth]{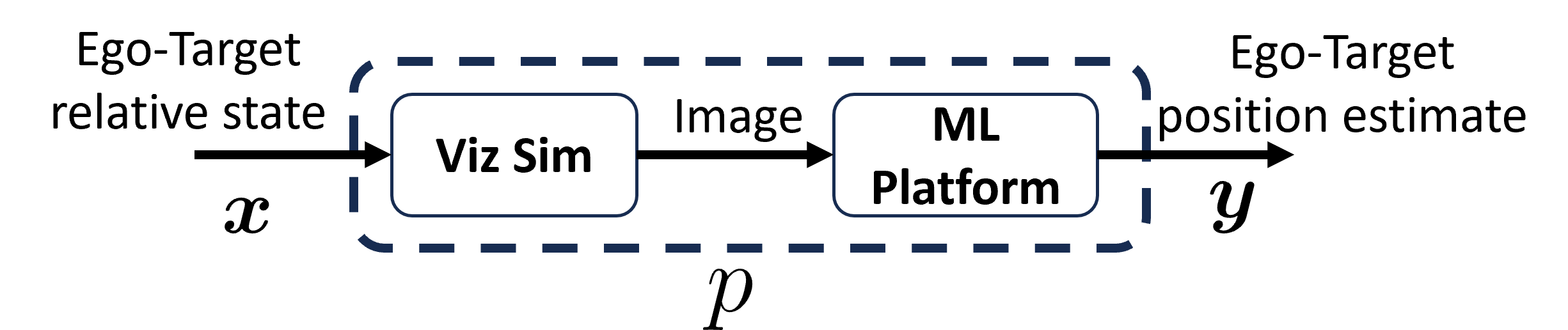}
    \caption{Schematic showing relative state to position estimate map}
    \label{fig:pmap}
\end{figure}

 The experiment design requires characterizing the uncertainty as a function of state-dependent parameters. We now need to create the state-to-output uncertainty map as discussed in Eq.~\ref{Eq:ellipse}, and approximate its upper bounding convex envelope. 
 In this setting, the perception map is denoted as $\bm y = p(\bm x)$, where $p$ denotes the process of generating the image corresponding to state $\bm x$ and performing inference on the image to generate an uncertain pose output $\bm y$.
 Fig.~\ref{fig:pmap} depicts this perception map, that operates on the agent's state $\bm x$. The Ego-Target pair is simulated and their relative state is defined by $\bm x$. An image of the Target is generated by process described in Fig.~\ref{fig:pmap}. This represents the sensors in any other application user wants to define. The ML estimator performs inference on this image to generate the output $\bm y = p(x) \sim U(\mathcal{Y}_{\bm x})$, with uniformly distributed uncertainty around $h(x)$ as the nominal output. The following assumption is now invoked:
 \begin{as}
     The perception map $\bm y = p(\bm x)$, samples the measurements with a uniform distribution from a set $\mathcal{Y_x}$. The size of the uncertainty set is given by a locally and globally bounded function $\Lambda(\mathcal{Y}_{\bm x})$.
 \end{as}
This is assumption implies that an outlier-rejection algorithm has been applied to the perception map output, keeping the output uncertainty bounded. Moreover, the map should also be immune to adversarial noise, which can affect the output when the image input is not generated strictly from the same distribution as the training dataset~\cite{adv}. For the case of satellite pose estimation scenario, the adversarial noise is non-existent as the simulator platform is used for testing and training  but can occur in practical setting.

For a state $\bm x$, the uncertainty is defined as $\bm z_{\bm x} = p(\bm x) - h(\bm x)$. Here, $h$ is the nominal output map for the system process; for the rendezvous for example, it is the relative position measurement for the state $\bm x$ without any uncertainty. With the bounded assumption, $\bm z_{\bm x}$ is said to be a uniformly distributed random variable. The 2-norm bound on $\bm z_x$ is the size of the uncertainty represented by the function $\Lambda(\mathcal{Y}_{\bm x})$, mentioned in Assumption~\ref{as:eigs}. The bounded random variable $\bm z_{\bm x}$ can be interpreted as being sub-Gaussian with  a variance $\sigma^2 = \Lambda(\mathcal{Y}_{\bm x})^2/4$~\cite{Vershynin_2018}. Concentration inequalities have been well established for sub-Gaussian random variables, with the consequence that their respective norms can be estimated for $N$ randomly drawn samples, as,
\begin{align}
    \bm d_x = \frac{1}{N}\sum\limits_{i=1}^N\lVert \bm z^i_{\bm x}\rVert, \label{eq:norm}
\end{align} with probability $1-\delta$, where $\delta \propto 1/N$~\cite{Vershynin_2018}. The approximation concentrates at $\sigma^2\sqrt{n}$, where $n$ is the dimension of $\bm z$. Therefore, the max bound can be computed as $\widehat{\Lambda}(\mathcal{Y}_{\bm x}) = \sqrt{\bm d_x/\sqrt{n}}$. The datapoints $z^i_{\bm x} = p(x) - h(x)$ are sampled by randomizing parameters of perception map $p$ that introduce the uncertainty.

 For the validation dataset, the enveloping function is generated by computing $\widehat{\Lambda}(\mathcal{Y}_{\bm x})$ for $M$ samples of $\bm{x} \in \mathcal{X}$. Here, $\bm x$ is sampled from $\mathcal{X}$ by generating a dense discretization on $\mathcal{X}$. For the dataset $(\bm x_i,\widehat{\Lambda}(\mathcal{Y}_{\bm x_i})), \; i=1,\ldots,M$, a tight quadratic envelope function of the form $\bm x^\top Q \bm x + b^\top \bm x + c$, where $Q\in \mathbb{S}^{n_x}_+, b\in \mathbb{R}^{n_x}, c\in \mathbb{R}$ is identified by solving the regression problem: 
 \begin{align}
     \min\limits_{Q\succeq0} &\quad \sum\limits_{i=1}^M\lVert \bm (x_i^\top Q\bm x_i + b^\top \bm x_i + c) - \widehat{\Lambda}(\mathcal{Y}_{\bm x_i}) \rVert_2^2 \label{eq:envqp} \\
     \mbox{s.t} &\quad \bm x_i^\top Q\bm x_i + b^\top \bm x_i + c - \widehat{\Lambda}(\mathcal{Y}_{\bm x_i}) \geq 0, \;\quad i = 1,\ldots,M.
 \end{align}
 
This problem can be compactly written as,
\begin{align}
     \min\limits_{q} &\quad \lVert L - Xq \rVert_2^2 \label{eq:envlp} \\
     \mbox{s.t} &\quad Xq - L \geq 0,\\
     &\quad q > 0.
 \end{align}
Note that the quadratic term has been expanded as $\bm x_i^\top Q\bm x_i = q_{11}x_1^2 + q_{12}x_1x_2 + \ldots+ q_{nn}x_n^2 + b_1 x_1 + \ldots b_n x_n + c := \bm q^\top \bm x_i'$, for $q = [q_{11}, q_{12},\ldots,q_{nn}, b_1,\ldots, b_n, c],\;  x_i' = [x_1^2,x_1x_2,\ldots,x_n^2, x_1,\ldots,x_n, 1 ]$, converting the quadratic term into its linear form for an $n$ dimensional $\bm x_i$. The data matrices are defined as $L = [\widehat{\Lambda}(\mathcal{Y}_{\bm x_1}),\ldots,\widehat{\Lambda}(\mathcal{Y}_{\bm x_M})]^\top$ and $X = [\bm x_1'^\top,\ldots,\bm x_M'^\top]^\top$. The constrained least-square problem can now be solved as an linear program, in order to generate the enveloping function, parametrized by $q$. 

\begin{figure}[t!]
\centering
\subcaptionbox{Representation of sampling strategy for validation step showing the Ego positions(blue) relative to Target(green).\label{fig:suna}}
{\includegraphics[width = 0.45\linewidth]{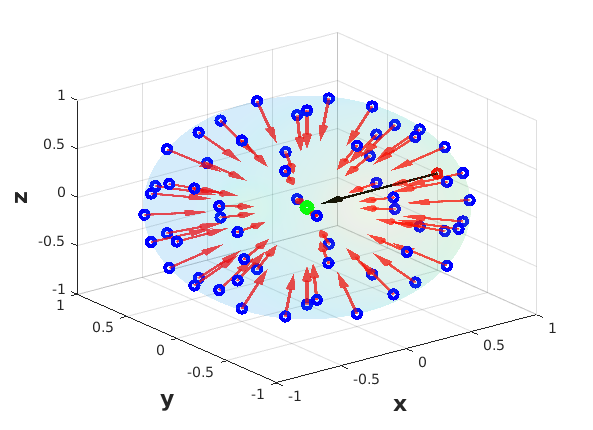}}\hfil
\subcaptionbox{ 2D representation of envelope(black) compared to maximum output uncertainty size vs Sun angle.\label{fig:sunb}}
{\includegraphics[width=0.45\linewidth]{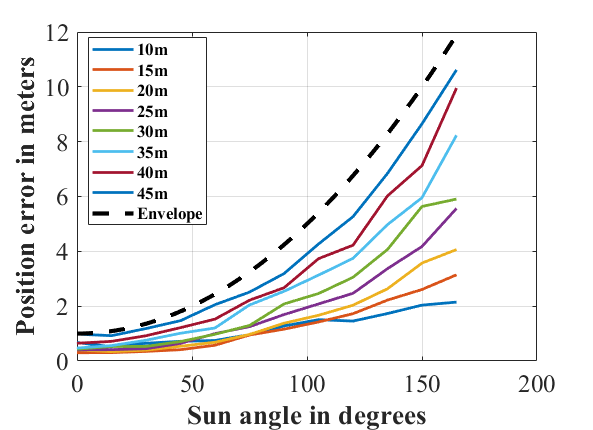}}
\caption{Validation experiment for fitting envelope function to output uncertainty.}\label{fig:sun}
\end{figure}

In this context, the system is given by Eq. (\ref{Eq:HCWx})-(\ref{Eq:HCWz}), for the state $\bm x$ representing the relative position and velocity of the Ego and the output $\bm y$ representing the relative position, as captured by the ML estimator. 
In order to generate a data-point of the dataset $(\bm x_i,\widehat{\Lambda}(\mathcal{Y}_{\bm x_i})$, a primary dataset is of $N$ points is collected, to solves for $\widehat{\Lambda}(\mathcal{Y}_{\bm x_i})$ using Eq. (\ref{eq:norm}) at each $\bm x_i$.  
The state $\bm x_i$ is sampled by densely discretizing the space $\mathcal{X}$ with a $n_x = 6$ dimensional grid.
For computing Eq. (\ref{eq:norm}) at each $\bm x_i$, we obtain a point $\bm z_{x_i}^j$ by running inference using the process shown in Fig.~\ref{fig:pmap}. The $N$ samples $\bm z_{x_i}^j$ are generated by varying the parameters for the simulation environment, that add randomness to image generated at $\bm x_i$. For each $\bm x_i$, the parameters that add uncertainty to the measurements have been classified as Target rotations, camera noise, and the background (Earth placement in the image), as detailed in~\cite{Ref:Becktor22}. A random Target rotation, Earth position and camera noise is sampled under feasible conditions to compute $\bm z_{\bm x_i}^j$. 
The process is repeated for $M$ points to generate the secondary dataset ($\widehat{\Lambda}(\mathcal{Y}_{\bm x_i}),\bm x_i$).
 Fig.~\ref{fig:suna} represents the points $\bm x_i$ sampled from $\mathcal{X}$ within the same distance from the Target (shown in green); this is shown as points (shown in blue) on a sphere. The Ego positions are shown in blue dots show the point $\bm x_i$ while the camera points towards the Target as shown by the red arrows. The Sun direction is shown by the black arrow.
 The secondary dataset $(\bm x_i,\widehat{\Lambda}(\mathcal{Y}_{\bm x_i})), \; i=1,\ldots,M$, is then populated with $\bm x_i$, sampled from $\mathcal{X}$ for feasible distances from the Target. The enveloping function of the quadratic form is now identified by solving Eq. (\ref{eq:envlp}). The solution to Eq. (\ref{eq:envlp}) gives the $n_{\bm x} = 6$ dimensional, quadratic enveloping function. Fig.~\ref{fig:sunb} represents a dimension of $\Lambda(\mathcal{Y}_{\bm x_i}) = \bm x^\top Q \bm x$ perpendicular to the Sun vector. The Sun angle is defined as the angle between the camera vector and the Sun vector. Fig.~\ref{fig:sunb} shows the representation of the quadratic envelope projected onto the 2D plane in-line with the Sun-angle. The dashed line covering the max error bounds represents the envelope. The solid lines represent variation of the max error at Ego positions with the same distance from the Target, with respect to their Sun-angles. Fig.~\ref{fig:sunb} also illustrates that uncertainty scales sharply as the distance to the Ego increases, as keypoint visibility reduces with distance. This is the main ``feature'' of this key-point based perception map that is exploited in our estimation aware planning.

Note that the data gathering step of the above process and the corresponding
uncertainty parametrization are influenced by the simulation platform and the problem scenario. For other applications where an estimation-aware trajectories need to be generate, the data can be collected in a similar manner by sampling the state space and the uncertainty domain. When the uncertainty set follows 
Assumption~\ref{app:sim}, the validation experiments can be used to define an enveloping function that upper bound the state-dependent uncertainty set.

To summarize, the enveloping function is generated by collecting a primary dataset that solves Eq. (\ref{eq:envlp}) to give $\widehat{\Lambda}(\mathcal{Y}_{\bm x})$ at each $\bm x$. Then  a secondary dataset is collected to compute $\widehat{\Lambda}(\mathcal{Y}_{\bm x_i})$ at each $\bm x_i$. Here $M$ such samples of $\bm x$ are collected from dense discretization of $\mathcal{X}$. This secondary dataset is used to solve Eq. (\ref{eq:envqp}) that defines the envelope. In principle any alternate method that computes a tight approximation described in Eq. (\ref{eq:lambdaapprox}) can be employed for the estimation-aware approach in this paper to be valid.
 

\bibliography{main}

\end{document}